%
%
%
%
%
\magnification=\magstep1
%
%
\newcount\parano
\newcount\eqnumbo
\newcount\thmno
\nopagenumbers
\binoppenalty=10000
\relpenalty=10000
\newcount\notitle
\notitle=1
%
%
\font\titlefont=cmssbx10 scaled \magstep1
\font\bigfont=cmr12
\font\eightrm=cmr8
\font\bbr=msbm10
\font\sbbr=msbm7
\font\mbbr=msbm9
\font\trm=cmr10 at 10truept
\font\tsl=cmsl10 at 10truept
\font\smallcaps=cmcsc10 at 10truept
\font\msam=msam10
%

%
%
\def\smallsect #1. #2\par{\bigbreak\noindent{\bf #1.}\enspace{\bf #2}\par
        \global\parano=#1\global\eqnumbo=1\global\thmno=1
        \nobreak\smallskip\nobreak\noindent\message{#2}}
\def\newthm#1 #2: #3{\xdef #2{\number\parano.\number\thmno}
        \global \advance \thmno by 1
        \medbreak\noindent
        {\bf #1 #2:}\if(#3\thmp\else\thmn#3\fi}
\def\thmp #1) { (#1)\thmn{}}
\def\thmn#1#2\par{\enspace{\it #1#2}\par
        \ifdim\lastskip<\medskipamount \removelastskip\penalty 55\medskip\fi}
\def\qedn{\thinspace\null\nobreak\hfill\hbox{\square}\par\medbreak}
\def\pf{\ifdim\lastskip<\smallskipamount \removelastskip\smallskip\fi
        \noindent{\it Proof\/}:\enspace}
\def\neweq#1$${\xdef #1{(\number\parano.\number\eqnumbo)}
        \eqno #1$$
        \global \advance \eqnumbo by 1}
\def\forevery#1#2$${\displaylines{\let\neweq=\newforclose
        \hfilneg\rlap{$\qqquad\forall#1$}\hfil#2\cr}$$}
\def\newforclose#1{
        \xdef #1{(\number\parano.\number\eqnumbo)}
        \hfil\llap{$#1$}\hfilneg
        \global \advance \eqnumbo by 1}
\def\itm#1{\par\indent\llap{\rm #1\enspace}\ignorespaces}
%
%
\gdef\begin #1\par{\xdef\titolo{#1}
\centerline{\titlefont\titolo}
\bigskip
\centerline{\bigfont by \autore}
\medskip
\centerline{\indirizzo}
\centerline{\email}
\centerline{\classificazione}
\medskip
\centerline{\data}
\bigskip\bigskip}
\headline={\ifodd\pageno\rhead\else\lhead\fi}
\def\rhead{\ifnum\pageno=\notitle\hfill\else\hfill\eightrm\titolo\hfill
\folio\fi}
\def\lhead{\ifnum\pageno=\notitle\hfill\else\eightrm\folio\hfill\autore\hfill
\fi}
%
%

\def\bar#1{\overline{#1}}
\let\de=\partial
\def\eps{\varepsilon}
\def\phe{\varphi}
\def\Hol{\mathop{\rm Hol}\nolimits}
\def\Re{\mathop{\rm Re}\nolimits}
\def\Im{\mathop{\rm Im}\nolimits}
\def\rKlim{\mathop{\hbox{$\tilde K$-$\lim$}}\limits}
\def\id{\mathop{\rm id}\nolimits}
\def\R{{\mathchoice{\hbox{\bbr R}}{\hbox{\bbr R}}{\hbox{\sbbr R}}
{\hbox{\sbbr R}}}}
\def\C{{\mathchoice{\hbox{\bbr C}}{\hbox{\bbr C}}{\hbox{\sbbr C}}
{\hbox{\sbbr C}}}}
\def\N{{\mathchoice{\hbox{\bbr N}}{\hbox{\bbr N}}{\hbox{\sbbr N}}
{\hbox{\sbbr N}}}}
\def\qqquad{\quad\qquad}
\def\interior #1{\mathaccent"7017 #1}
\def\mezzo{{\textstyle{1\over2}}}
\def\n{\hbox{$|\mkern-2mu|\mkern-2mu|$}}
\def\ln{\mathopen{\n}}
\def\rn{\mathclose{\n}}
\def\square{{\msam\char"03}}
\def\void{\hbox{\bbr\char"3F}}
%
%
\newbox\bibliobox
\def\setref #1{\setbox\bibliobox=\hbox{[#1]\enspace}
        \parindent=\wd\bibliobox}
\def\biblap#1{\noindent\hang\rlap{[#1]\enspace}\indent\ignorespaces}
\def\art#1 #2: #3! #4! #5 #6 #7-#8 \par{\biblap{#1}#2: #3.
        #4~{\bf #5} \hbox{#7--#8} (#6)\par\smallskip}
\def\book#1 #2: #3! #4 \par{\biblap{#1}#2: #3. #4\par\smallskip}
%
%
\def\autore{Marco Abate}
\def\indirizzo{\trm Dipartimento di Matematica, Universit\`a di Ancona, via
Brecce Bianche, 60131 Ancona, Italy}
\def\email{\trm Fax: 39/71/220.4870. E-mail: abate@anvax1.cineca.it}
\def\classificazione{{\tsl Mathematics Subject Classification (1991):} \trm
32A40, 32H40,  32H15}
\def\data{December 1996}
%
%
\begin The Julia-Wolff-Carath\'eodory theorem in polydisks

\centerline{\smallcaps Abstract}
\smallskip
{\narrower\narrower\normalbaselineskip=11pt\normalbaselines

  \trm\noindent The classical Julia-Wolff-Carath\'eodory theorem
  gives a condition ensuring the existence of the non-tangential
  limit of both a bounded holomorphic function and its derivative
  at a given boundary point of the unit disk in the complex
  plane. This theorem has been generalized by Rudin to
  holomorphic maps between unit balls in~$\hbox{\mbbr C}^n$, and
  by the author to holomorphic maps between strongly
  (pseudo)convex domains. Here we describe
  Julia-Wolff-Carath\'eodory theorems for holomorphic maps
  defined in a polydisk and with image either in the unit disk,
  or in another polydisk, or in a strongly convex domain. One of
  main tool for the proof is a general version of Lindel\"of's
  principle valid for not necessarily bounded holomorphic
  functions.

}

\smallsect 0. Introduction

The classical Fatou theorem says that a bounded holomorphic
function~$f$ defined on the unit disk~$\Delta\subset\C$ admits
non-tangential limit at almost every point of~$\de\Delta$, but it
does not say anything about the behavior of~$f(\zeta)$ as~$\zeta$
approaches a specific point~$\sigma$ of the boundary.  Of course,
to be able to say something in this case one needs some
hypotheses on~$f$. For instance, one can assume that, in a very
weak sense, $f(\zeta)$ goes to the boundary of the image of~$f$
as~$\zeta$ goes to~$\sigma$. This leads to the classical {\sl
  Julia's lemma:}

\newthm Theorem \zJulia: (Julia [Ju1]) Let $f\colon\Delta\to\Delta$ be a
bounded holomorphic function such that
$$\liminf_{\zeta\to\sigma}{1-|f(\zeta)|\over1-|\zeta|}=\alpha<+\infty
        \neweq\eqzalp$$
for some $\sigma\in\de\Delta$. Then $f$ has non-tangential limit $\tau\in\de
\Delta$ at~$\sigma$, and furthermore
$${|\tau-f(\zeta)|^2\over1-|f(\zeta)|^2}\le\alpha\,{|\sigma-\zeta|^2\over1-
        |\zeta|^2}\neweq\eqzJ$$
for all $\zeta\in\Delta$.

This statement has a very interesting geometrical interpretation. The {\sl
horocycle} $E(\sigma,R)\subset\Delta$ of {\sl center}~$\sigma\in\de\Delta$ and
{\sl radius}~$R>0$ is the set
$$E(\sigma,R)=\left\{\zeta\in\Delta\biggm|{|\sigma-\zeta|^2\over1-|\zeta|^2}
        <R\right\}.$$
Geometrically, $E(\sigma,R)$ is an euclidean disk of radius $R/(R+1)$
internally tangent to~$\de\Delta$ at~$\sigma$. Therefore~\eqzJ\ becomes
$f\bigl(E(\sigma,R)\bigr)\subseteq E(\tau,\alpha R)$ for all $R>0$, and the
existence of the non-tangential limit more or less follows from~\eqzJ\ and from
the fact that horocycles touch the boundary in exactly one point.

A horocycle can be thought of as the limit of Poincar\'e disks of
fixed euclidean radius and centers going to the boundary; in a
sense, thus, we can interpret horocycles as Poincar\'e disks
centered at the boundary, and Julia's lemma as a Schwarz-Pick
lemma at the boundary. This suggests that $\alpha$ might be
considered as a sort of dilatation coefficient: $f$ expands
horocycles by a ratio of~$\alpha$. If $\sigma$ were an internal
point and $E(\sigma,R)$ an infinitesimal euclidean disk actually
centered at~$\sigma$, one would be tempted to say that~$\alpha$
is the absolute value of the derivative of~$f$ at~$\sigma$.  This
is exactly the content of the classical {\sl
  Julia-Wolff-Carath\'eodory theorem:}

\newthm Theorem \zJWC: Let $f\colon\Delta\to\Delta$ be a
bounded holomorphic function such that
$$\liminf_{\zeta\to\sigma}{1-|f(\zeta)|\over1-|\zeta|}=\alpha<+\infty$$
for some $\sigma\in\de\Delta$, and let $\tau\in\de\Delta$ be the non-tangential
limit of~$f$ at~$\sigma$. Then both the incremental ratio $\bigl(\tau-f(\zeta)
\bigr)\big/(\sigma-\zeta)$ and the derivative~$f'(\zeta)$ have
non-tangential limit~$\alpha\bar\sigma\tau$ at~$\sigma$.

So condition \eqzalp\ forces the existence of the non-tangential limit of
both~$f$ and its derivative at~$\sigma$. This theorem results from the work of
several people: Julia~[Ju2], Wolff~[Wo], Carath\'eodory~[C], Landau and
Valiron~[L-V], R. Nevanlinna~[N] and others (see Burckel~[B] and~[A4] for
proofs,
history and applications).

Before describing what happens for
functions of several complex variables, it is worthwhile to take a closer look
to the notion of non-tangential limit in~$\Delta$. First of all, the
non-tangential limit can be defined in two equivalent ways. We can say that a
function $f\colon\Delta\to\C$ has {\sl non-tangential limit}~$L\in\C$
at~$\sigma\in\de\Delta$ if~$f\bigl(\gamma(t)\bigr)\to L$ as~$t\to 1^-$ for
every curve $\gamma\colon[0,1)\to\Delta$ such that $\gamma(t)\to\sigma$
non-tangentially as~$t\to1^-$. Or, we can say that $f$ has non-tangential
limit~$L\in\C$ if $f(\zeta)\to L$ as~$\zeta\to\sigma$ staying inside any {\sl
Stolz region}~$H(\sigma, M)$ of {\sl vertex}~$\sigma$ and~{\sl
amplitude}~$M>1$, where
$$H(\sigma,M)=\left\{\zeta\in\Delta\biggm|{|\sigma-\zeta|\over1-|\zeta|}<M.
        \right\}.$$
Since Stolz regions are angle-shaped nearby the vertex~$\sigma$, and the angle
is going to~$\pi$ as~$M\to+\infty$, it is clear that these two definitions are
equivalent; as we shall see later on, this is not anymore the case in~$\C^n$
with~$n>1$.

The second observation is that to check the existence of a non-tangential limit
it suffices to study the function along one single curve. This is the content
of {\sl Lindel\"of's principle} (see Rudin~[R] for a modern proof):

\vbox{\newthm Theorem \zL: (Lindel\"of [Li]) Let $f\colon\Delta\to\C$ be a
bounded holomorphic function. Assume there are $\sigma\in\de\Delta$ and a
continuous curve~$\gamma\colon[0,1)\to\Delta$ with $\gamma(t)\to\sigma$
as~$t\to1^-$ so that
$$\lim_{t\to1^-}f\bigl(\gamma(t)\bigr)=L\in\C.$$
Then $f$ has non-tangential limit~$L$ at~$\sigma$.\par}

All these results have been generalized in several ways to functions and maps
defined in domains of~$\C^n$ (see the introductions to~[A3,~5] for a brief
historical summary); it turns out that the best way to express the
generalizations (or, at least, the best way according to the point of view
discussed in this paper) is in terms of geometric function theory.

The aim of geometric function theory is to describe how the geometrical shape
of a domain influences the analytical behavior of holomorphic functions defined
on it. A typical (and of the utmost importance) result of this kind is the
Schwarz-Pick lemma: a holomorphic self-map of~$\Delta$ is necessarily a
contraction with respect to the Poincar\'e distance (which is a
geometrical object depending only on the shape of~$\Delta$).

The Kobayashi distance of a domain~$D\subset\subset\C^n$ is a
direct generalization of the Poincar\'e distance of~$\Delta$
(see, e.g., [K1,~2],~[J-P] and~[A4] for definition and
properties); in particular, it is contracted by holomorphic maps.
The main claim of this paper is that all the results discussed
before (Julia's lemma, Julia-Wolff-Carath\'eodory theorem,
Lindel\"of's principle) are better understood if expressed using
the Kobayashi distance (and its infinitesimal analogue, the
Kobayashi metric). In previous papers (see~[A3,~5]) we have
already shown how to do so in strongly convex and strongly
pseudoconvex domains; here we shall show how the same techniques
can be applied to polydisks (and conceivably to other domains as
well).

The first step consists in understanding the correct
generalization to several variables of the notion of
non-tangential limit. As first discovered by Kor\'anyi and Stein
discussing Fatou's theorem in~$\C^n$ ([Ko, K-S, S]), the natural
approach regions for the study of boundary behavior of
holomorphic functions are not cones, but regions approaching the
boundary non-tangentially along the normal direction, and
tangentially (at least parabolically) along the complex
tangential directions. Kor\'anyi and Stein defined their approach
regions in euclidean (local) terms; later, Krantz~[Kr] suggested
a definition using both the Kobayashi distance and the euclidean
normal line. In~[A3] we generalized Stolz regions as follows:
if~$D\subset \subset\C^n$ is a bounded domain and~$z_0\in D$,
then the ({\sl small\/}) {\sl Kor\'anyi region}~$H_{z_0}(x,M)$ of
{\sl vertex}~$x\in\de D$, {\sl pole}~$z_0\in D$ and {\sl
  amplitude}~$M>1$ is given by
$$H_{z_0}(x,M)=\bigl\{z\in D\bigm|\limsup_{w\to
  x}[k_D(z,w)-k_D(z_0,w)]+k_D(z_0,z)<\log M\bigr\},$$ where $k_D$
is the Kobayashi distance of~$D$. The rationale behind this
definition is the following: the $\limsup$ measure the
``distance'' (normalized as to be always finite, possibly
negative) of~$z$ from the boundary point~$x$; so~$z\in
H_{z_0}(x,M)$ if and only if the average between the ``distance''
of~$z$ from~$x$ and the distance of~$z$ from the pole~$z_0$ is
bounded by~$\mezzo\log M$. By the way, changing the pole amounts
to shifting the amplitudes, and thus it is completely irrelevant.

It is easy to check that the Stolz
region~$H(\sigma,M)\subset\Delta$ is exactly the Kor\'anyi region
of vertex~$\sigma$, pole the origin and amplitude~$M$.  Therefore
we shall say that a function~$f\colon D\to\C$ has {\sl
  $K$-limit}~$L\in\C$ at~$x\in\de D$ if~$f(z)\to L$ as~$z\to x$
inside any Kor\'anyi region of vertex~$x$.

Since if $D$ is strongly pseudoconvex ([A5]) our Kor\'anyi regions are
comparable (actually equal if $D=B^n$, the unit ball of~$\C^n$) with
Kor\'anyi-Stein's and Krantz's approach regions (with the latter even when
$D$ is
a polydisk; see Remark~1.5), our $K$-limit coincides with their admissible
limit, but our approach regions are completely described in terms of the
Kobayashi distance only.

We remarked before that the non-tangential limit in~$\Delta$ can be defined in
two different ways. We just generalized the first way; \v Cirka~[\v C]
discovered that the second one leads to a different notion of limit, which is
the one needed to generalize Lindel\"of's principle.

Let $D\subset\subset\C^n$ be a bounded domain, and $x\in\de D$. A
{\sl $x$-curve} is a continuous curve~$\sigma\colon[0,1)\to D$
such that $\sigma(t)\to x$ as~$t\to1^-$. Then \v Cirka showed
that the correct version of Lindel\"of's principle should be
something like this: if $f\colon D\to\C$ is a bounded holomorphic
function such that $f\bigl(\sigma^o(t)\bigr)\to L\in\C$
as~$t\to1^-$ for a given $x$-curve~$\sigma^o$ belonging to a
specific class of $x$-curves, then $f\bigl(\sigma(t)\bigr)\to L$
as~$t\to1^-$ for all $x$-curves~$\sigma$ belonging to a (possibly
sub)class of $x$-curves.

The main point here clearly is the identification of the correct class of
curves. Several authors (see, e.g., [\v C, C-K, D, D-Z, Kh]) have suggested
several possibilities, more or less intrinsic. In [A3] we found a general
technique (inspired by~[C-K]) to generate such classes. Let
$D\subset\subset\C^n$
be a bounded domain, and $x\in\de D$. A {\sl projection device} at~$x$ is given
by the following data: a holomorphic immersion~$\phe_x\colon\Delta\to D$
extending continuously to the boundary so that~$\phe_x(1)=x$; a
neighborhood~$U_0$ of~$x$ in~$\C^n$; and a device associating to every
$x$-curve~$\sigma$ contained in~$U_0$ a
$x$-curve~$\sigma_x$ contained in~$U_0\cap\phe_x(\Delta)$. For instance, if~$D$
is strongly convex containing the origin, we can take as $\phe_x$ a
parametrization of the intersection of~$D$ with the complex line~$\C x$ passing
through~$0$ and~$x$, and we can associate to every~$x$-curve its orthogonal
projection into~$\C x$. In~[A5] we described two other projection devices,
another one is introduced in this paper (see Section~1), and many others are
easily devised.

Given a projection device at~$x\in\de D$, we say that a $x$-curve
$\sigma$ is {\sl restricted} if the
curve~$\tilde\sigma_x=\phe^{-1}_x\circ\sigma_x$ goes to~$1$
non-tangentially in~$\Delta$, and that the $x$-curve~$\sigma$ is
{\sl special} if~$k_D\bigl(\sigma(t),\sigma_x(t)\bigr)\to0$
as~$t\to1^-$. We say that a function~$f\colon D\to\C$ has {\sl
  restricted} \hbox{\sl $K$-limit}~$L$ at~$x$ (with respect to
the given projection device) if $f\bigl(\sigma(t)\bigr)\to L$
as~$t\to1^-$ along any special restricted $x$-curve~$\sigma$.
Then

\newthm Theorem \zLA: ([A3]) Let $D\subset\subset\C^n$ be a bounded domain
equipped with a projection device at~$x\in\de D$. Let $f\colon D\to\C$ be a
bounded holomorphic function, and assume there is a special
$x$-curve~$\sigma^o$ such that
$$\lim_{t\to1^-}f\bigl(\sigma^o(t)\bigr)=L\in\C.$$
Then $f$ has restricted $K$-limit~$L$ at~$x$.

For instance, if we use the projection device defined before, Theorem~\zLA\
recovers \v Cirka's results for the case of strongly convex domains.

We should remark that whereas Theorem~\zLA\ holds in this very general setting,
to generalize Theorem~\zJWC\ we shall need a similar result for functions which
are only bounded in Kor\'anyi regions. At present there is no general proof of
such a statement, but only (usually hard; cf. Theorem~2.2) proofs working
slightly differently in each case, depending on the actual geometry of the
domain under consideration. It would be interesting to have a more general
proof.

The next step is the generalization of Julia's lemma. Recalling the description
of the horocycles as Poincar\'e disks at the boundary, it is natural to define
in a domain $D\subset\subset\C^n$ the {\sl small horosphere}~$E_{z_0}(x,R)$ of
{\sl center}~$x\in\de D$, {\sl pole}~$z_0\in D$ and {\sl radius}~$R>0$ by
$$E_{z_0}(x,R)=\bigl\{z\in D\bigm|\limsup_{w\to x}[k_D(z,w)-k_D(z_0,w)]<
        \mezzo\log R\bigr\};$$
the {\sl big horosphere}~$F_{z_0}(x,R)$ is analogously defined replacing the
$\limsup$ by a $\liminf$.

In condition \eqzalp, both numerator and denominator can be interpreted as
distances from the boundary. Since, when $D$ is complete hyperbolic and~$z_0\in
D$ is chosen once for all, $\exp\bigl(-2k_D(z_0,z)\bigr)$ is going to zero
exactly when~$z$ is going to the boundary (and it behaves exactly as the
euclidean distance from the boundary when $D$ is strongly pseudoconvex;
see~[A1]), we can use it as a replacement for the euclidean distance from the
boundary. This leads to the following version of Julia's lemma:

\newthm Theorem \zJA: ([A3]) Let $D_1\subset\subset\C^n$,
$D_2\subset\subset\C^m$
be two bounded domains, with $D_1$ complete hyperbolic. Fix~$z_1\in D_1$
and~$z_2\in D_2$. Let $f\colon D_1\to D_2$ be a holomorphic map such that
$$\liminf_{w\to
  x}\bigl[k_{D_1}(z_1,w)-k_{D_2}\bigl(z_2,f(w)\bigr)\bigr]=
\mezzo\log\alpha<+\infty\neweq\eqzaJ$$ for some $x\in\de D$ and
$\alpha>0$. Then there exists a $y\in\de D_2$ such that
$$\forevery{R>0}f\bigl(E_{z_1}(x,R)\bigr)\subseteq F_{z_2}(y,\alpha R).$$
Furthermore, if $\bar{F_{z_2}(y,R)}\cap\de D_2=\{y\}$ for all~$R>0$ (for
instance, if $D_2$ is strongly pseudoconvex) then $f$ has $K$-limit~$y$ at~$x$.

It should be remarked that in [A3] this theorem was stated only for
self-maps of a complete hyperbolic domain, but the same proof goes through in
the general case too.

Very often, \eqzaJ\ is equivalent to
$$\liminf_{w\to x}{d\bigl(f(w),\de D_2\bigr)\over d(w,\de D_1)}<+\infty,$$
and so this result encompasses most of the other known generalizations of
Julia's
lemma (see, e.g., [M, H, W]).

We are left with the generalization of the Julia-Wolff-Carath\'eodory theorem.
With respect to the one-dimensional case there is an obvious difference:
instead of only one derivative we have to consider a whole (Jacobian) matrix of
them, and there is no reason they should all behave in the same way. And
indeed, as first showed by Rudin~[R] in~$B^n$, they do not. It turns out that
the geometry of the domain (and, in particular, the boundary behavior of the
Kobayashi metric) plays a main role in determining the correct statement of a
Julia-Wolff-Carath\'eodory theorem.

Rudin proved the following Julia-Wolff-Carath\'eodory theorem for the unit ball
of~$\C^n$:

\newthm Theorem \zJWCR: (Rudin [R]) Let $f\colon B^n\to B^m$ be a holomorphic
map such that
$$\mezzo\log\liminf_{w\to
  x}{1-\|f(w)\|\over1-\|w\|}=\liminf_{w\to x}
\bigl[k_{B^n}(0,w)-k_{B^m}\bigl(0,f(w)\bigr)\bigr]=\mezzo\log\alpha<+\infty,$$
for some $x\in B^n$, where $\|\cdot\|$ denote the euclidean norm.
Let $y\in\de B^m$ be the $K$-limit of~$f$ at~$x$, as given by
Theorem~\zJA. Equip~$B^n$ with the projection device defined
before in terms of the orthogonal projection, set
$f_y(z)=\bigl(f(z),y \bigr)y$ (where $(\cdot\,,\cdot)$ is the
canonical hermitian product) and denote by~$df_z$ the
differential of~$f$ at~$z$. Then: {\smallskip \itm{(i)} The
  incremental ratio
  $\bigl[1-\bigl(f(z),y\bigr)\bigr]\big/[1-(z,x)]$ has restricted
  $K$-limit~$\alpha$ at~$x$; \itm{(ii)} The map
  $[f(z)-f_y(z)]/[1-(z,x)]^{1/2}$ has restricted $K$-limit~$0$
  at~$x$; \itm{(iii)} The function $\bigl(df_z(x),y\bigr)$ has
  restricted $K$-limit~$\alpha$ at~$x$; \itm{(iv)} The map
  $[1-(z,x)]^{1/2}d(f-f_y)_z(x)$ has restricted $K$-limit~$0$
  at~$x$;
\item{\rm(v)} If $x^\perp$ is any vector orthogonal to~$x$, the function
$\bigl(df_z (x^\perp),y\bigr)\big/[1-(z,x)]^{1/2}$ has restricted $K$-limit~$0$
at~$x$;
\item{\rm(vi)} If $x^\perp$ is any vector orthogonal to~$x$, the
  map $d(f-f_y)_z (x^\perp)$ is bounded in every\break\indent
  Kor\'anyi region of vertex~$x$.}

In~[A3,~5] we generalized this theorem to strongly convex and strongly
pseudoconvex domains, showing explicitly that the different behavior of the
partial derivatives is due to the different behavior of the Kobayashi metric
along normal directions or complex tangential directions.

Recently, Jafari~[J] made a preliminary study of the case of the polydisk
along the \v Silov boundary; unfortunately, one of his statements is wrong (see
Remarks~4.2 and~4.5). This paper is devoted to a full description
of what happens in the polydisk. Our main result (Theorem~4.1) deals with
bounded
holomorphic functions:

\newthm Theorem \zJWCA: Let $f\colon\Delta^n\to\Delta$ be a
holomorphic function, and~$x=(x_1,\ldots,x_n)\in\de\Delta^n$. Assume there
is~$\alpha>0$ such that
$$\mezzo\log\liminf_{w\to x}{1-|f(w)|\over1-\ln w\rn}=\liminf_{w\to
        x}\bigl[k_{\Delta^n}(0,w)-k_{\Delta}\bigl(0,f(w)\bigr)
        \bigr]=\mezzo\log\alpha<+\infty,$$
where $\ln\cdot\rn$ denote the sup-norm.
Let $\tau\in\de\Delta$ be the $K$-limit of~$f$ at~$x$, as given by
Theorem~\zJA. Equip $\Delta^n$ with the canonical projection device at~$x$ (to
be described in Section~1). Then:
{\smallskip
\item{\rm(i)}the incremental ratio $\bigl(\tau-f(z)\bigr)\big/\bigl(1-\tilde
p_x(z)\bigr)$ has restricted $K$-limit~$\alpha\tau$ at~$x$, where the
function~$\tilde p_x\colon\Delta^n\to D$ is defined in~$(1.2)$;
\item{\rm(ii)}If $|x_j|=1$, then the incremental ratio
$\bigl(\tau-f(z)\bigr)\big/(x_j-z_j)$ has restricted
$K$-limit~$\alpha\tau\bar{x_j}$ at~$x$;
\itm{(iii)} the partial derivative $\de f/\de x=df_z(x)$ has restricted
$K$-limit~$\alpha\tau$ at~$x$;
\itm{(iv)}If $|x_j|<1$ then the partial derivative $\de f/\de z_j$ has
restricted $K$-limit~$0$ at~$x$;
\itm{(v)}If $|x_j|=1$ then $\de f/\de z_j$ has restricted $K$-limit at~$x$.}

This theorem is proved in Section~4. Section~1 is devoted to define
Kor\'anyi regions and the canonical projection device in~$\Delta^n$. In
Section~2
we prove the generalization of Lindel\"of's principle, and in Section~3 of
Julia's lemma, along the lines described before. Finally, in Section~5 we shall
discuss some more Julia-Wolff-Carath\'eodory theorems for maps from~$\Delta^n$
into another polydisk, or into a strongly convex domains, showing how it is
possible to apply in several different
situations the machinery we developed here.

\smallsect 1. The canonical projection device and Kor\'anyi regions

As described in the introduction, to prove a Julia-Wolff-Carath\'eodory theorem
two tools are needed: a Lindel\"of principle, and a Julia's lemma. The aim of
this section is to introduce concepts and definitions needed by both of them.

We shall denote by $\|\cdot\|$ the euclidean norm on $\C^n$, and by $(\cdot\,,
\cdot)$ the canonical hermitian product
$$(z,w)=z_1\bar{w_1}+\cdots+z_n\bar{w_n}.$$
Let $\ln\cdot\rn$ denote the sup norm
$$\ln z\rn=\max_j\{|z_j|\};$$
the {\sl unit polydisk} $\Delta^n\subset\C^n$ is the unit ball for this norm,
that is $\Delta^n=\{z\in\C^n\mid\ln z\rn<1\}$.

The Kobayashi distance $k_{\Delta^n}$ of $\Delta_n$ is given by (see, e.g.,
[J-P, Example~3.1.8])
$$k_{\Delta^n}(z,w)=\mezzo\log{1+\ln\gamma_z(w)\rn\over1-\ln\gamma_z(w)\rn},$$
where $\gamma_z\colon\Delta^n\to\Delta^n$ is the automorphism of $\Delta^n$
given by
$$\gamma_z(w)
=\left({w_1-z_1\over1-\bar{z_1}w_1},\cdots,{w_n-z_n\over1-\bar{z_n}
    w_n}\right).$$ Since $t\mapsto\mezzo\log(1+t)/(1-t)$ is an
increasing function, we get
$$k_{\Delta_n}(z,w)
=\max_j\left\{\mezzo\log{1+\left|{w_j-z_j\over1-\bar{z_j}w_j}
    \right|\over1-\left|{w_j-z_j\over1-\bar{z_j}w_j}\right|}\right\}=
\max_j\,\{\omega(z_j,w_j)\},\neweq\eqkd$$ where $\omega$ is the
Poincar\'e distance on the unit disk $\Delta\subset\C$.

A {\sl complex geodesic} in a domain $D\subset\C^n$ is a a
holomorphic map $\phe\colon\Delta\to D$ which is an isometry with
respect to the Poincar\'e distance on $\Delta$ and the Kobayashi
distance on~$D$. It is easy to prove (see, e.g., [A4,
Proposition~2.6.10]) that $\phe\in\Hol(\Delta,\Delta^n)$ is a
complex geodesic iff at least one component of~$\phe$ is an
automorphism of~$\Delta$. In particular, given two points
of~$\Delta^n$ there is always a complex geodesic containing both
of them in its image, and given a point of~$\Delta^n$ and a
tangent vector there is always a complex geodesic containing the
point and tangent to the vector.

We shall need a particular class of complex geodesics in $\Delta^n$. Given
$x\in\de\Delta^n$, the {\sl complex geodesic associated} to~$x$ is the map
$\phe_x\in\Hol(\Delta,\Delta^n)$ given by
$$\phe_x(\zeta)=\zeta x.$$
Since there is at least one component~$x_j$ of $x$ with $|x_j|=1$,
every~$\phe_x$ is a complex geodesic, and we have a canonical family of complex
geodesics whose images fill up all of~$\Delta^n$ and meet only at the origin.

To better describe the boundary behavior of functions and maps we
shall need some more terminology. Let
$x=(x_1,\ldots,x_n)\in\de\Delta^n$; the {\sl \v Silov
  degree}~$d_x$ of~$x$ is the number of components of~$x$ with
absolute value~1:
$$d_x=\#\{j\mid|x_j|=1\}.$$
In particular, $d_x=n$ iff $x$ belongs to the {\sl \v Silov boundary} $(\de
\Delta)^n=\de\Delta\times\cdots\times\de\Delta\subset\de\Delta^n$
of~$\Delta^n$. A {\sl \v Silov component} of~$x$ is a component~$x_j$ such
that~$|x_j|=1$; an {\sl internal component} of $x$ is a component~$x_j$ such
that~$|x_j|<1$. The {\sl \v Silov part} $\check x=(\check x_1,\ldots,\check
x_n)\in\de\Delta^n$ of~$x$ is defined by
$$\check x_j=\cases{x_j&if $|x_j|=1$,\cr 0&if $|x_j|<1$;\cr}$$
notice that $(x,\check x)=d_x$. The {\sl internal part}
$\interior x\in\Delta^n$ of~$x$ is defined by~$\interior
x=x-\check x$, so that~$x=\check x+\interior x$. Finally, we
shall say that a tangent vector~$v=(v_1,\ldots,v_n)\in\C^n$ has
{\sl no \v Silov components with respect to}~$x$ if $v_j=0$
whenever~$|x_j|=1$.

We have already associated to each $x\in\de\Delta^n$ a complex
geodesic~$\phe_x$. The {\sl left-inverse function associated
  to}~$x$ is the holomorphic function $\tilde
p_x\colon\Delta^n\to\Delta$ given by
$$\tilde p_x(z)={1\over d_x}(z,\check x);\neweq\eqptx$$
the {\sl holomorphic retraction} associated to $x$ is the holomorphic
map~$p_x\colon\Delta^n\to\Delta^n$ given by
$$p_x(x)={1\over d_x}(z,\check x)x=\phe_x\circ\tilde p_x(z).$$
Notice that $\tilde p_x\circ\phe_x=\id_\Delta$ (hence the name
left-inverse), $p_x\circ p_x=p_x$ (hence the name retraction),
and $p_x\circ\phe_x=\phe_x$, so that $p_x$ is a retraction
of~$\Delta^n$ onto the image of~$\phe_x$.

A {\sl $x$-curve} is a continuous curve $\sigma\colon[0,1)\to\Delta^n$ such
that $\sigma(t)\to x$ as~$t\to 1^-$. We shall say that something happens {\sl
eventually} to a $x$-curve $\sigma$ if it happens to~$\sigma(t)$ for $t$ close
enough to~1. If $\sigma=(\sigma_1,\ldots,\sigma_n)$, then an {\sl internal
component} (respectively, a {\sl \v Silov component\/}) of~$\sigma$ is a
component~$\sigma_j$ such that $x_j$ is an internal (respectively, \v Silov)
component of~$x$.

We are now ready to define the projection device (see the
introduction) we shall use.  Let~$x\in\de\Delta^n$; the {\sl
  canonical projection device} at~$x$ is given by the complex
geodesic~$\phe_x$, the neighborhood~$U_0=\C^n$ and by associating
to any $x$-curve~$\sigma$ the $x$-curve
$$\sigma_x=p_x\circ\sigma={1\over d_x}(\sigma,\check x)x;$$ it is
a $x$-curve whose image is contained in~$\phe_x(\Delta)$, and it
does not depend on the internal components of~$\sigma$.
Furthermore,
$$\ln\sigma_x\rn={1\over d_x}|(\sigma,\check x)|,\neweq\eqnsx$$
and
$$\ln\sigma-\sigma_x\rn
=\max_j\left\{\left|\sigma_j-{(\sigma,\check x)\over d_x}
        x_j\right|\right\}.\neweq\eqnssx$$
We shall also set $\tilde\sigma_x=\tilde p_x\circ\sigma$, which is a 1-curve
in~$\Delta$ such that~$\sigma_x=\phe_x\circ\tilde\sigma_x$. Notice that~$\ln
\sigma_x\rn=|\tilde\sigma_x|$.

We shall say that a $x$-curve $\sigma$ is {\sl special} if
$$\lim_{t\to 1^-}k_{\Delta^n}\bigl(\sigma(t),\sigma_x(t)\bigr)=0.$$
This is equivalent to requiring that
$$\max_j\left\{\left|{\sigma_j-{1\over d_x}(\sigma,\check x)x_j\over
   1-\bar{\sigma_j}{1\over d_x}(\sigma,\check x)x_j}\right|\right\}\to 0.$$
Let $x_j$ be an internal component of~$x$. Then
$$\bar{\sigma_j}{1\over d_x}(\sigma,\check x)x_j\to|x_j|^2<1;$$
therefore $\sigma$ is special iff
$$\max_{|x_j|=1}\left\{\left|{\sigma_j-{1\over d_x}(\sigma,\check x)x_j\over
        1-\bar{\sigma_j}{1\over d_x}(\sigma,\check x)x_j}\right|\right\}\to 0.
        \neweq\eqsp$$
In particular, being special imposes no restrictions on the internal components
of a $x$-curve.

There is a geometric interpretation of special curves: a $x$-curve is special
if the \v Silov components approach the image of~$\phe_x$ faster than
the rate of approach of the projection~$\sigma_x$ to the boundary
of~$\Delta^n$. More precisely,

\newthm Proposition \Jafari: Let $x\in\de\Delta^n$. Then a $x$-curve $\sigma$
is special if and only if
$$\lim_{t\to 1^-}\max_{|x_j|=1}\left\{{|\sigma_j(t)-(\sigma_x)_j(t)|\over1-\ln
        \sigma_x(t)\rn}\right\}=0.\neweq\eqJs$$

\pf Assume first that \eqJs\ holds. Then for $j=1,\ldots,n$ we have
$$\left|1-\bar{\sigma_j}{1\over d_x}(\sigma,\check x)x_j\right|\ge 1-{1\over
        d_x}|(\sigma,\check x)|,$$
because $|\sigma_j|<1$ and $|x_j|\le 1$. Therefore
$$\max_{|x_j|=1}\left|{\sigma_j-{1\over d_x}(\sigma,\check x)x_j\over
  1-\bar{\sigma_j}{1\over d_x}(\sigma,\check x)x_j}\right|\le\max_{|x_j|=1}
        \left\{{|\sigma_j-(\sigma_x)_j|\over1-\ln\sigma_x\rn}\right\},$$
by \eqnsx\ and \eqnssx, and so $\sigma$ is special, by \eqsp.

Conversely, assume \eqJs\ is not satisfied. Write
$\sigma=\sigma_x+\alpha$, where $\alpha\colon[0,1)\to\C^n$ is such
that~$(\alpha,\check x)\equiv 0$ and $\alpha(t)\to0$ as~$t\to1^-$. Since
\eqJs\ is not satisfied, there are a \v Silov component $\sigma_j$, an
$\eps>0$ and a sequence $t_k\to1^-$ such that
$${|\alpha_j(t_k)|\over1-|\tilde\sigma_x(t_k)|^2}
\ge{|\sigma_j(t_k)-(\sigma_x)_j
        (t_k)|\over2(1-\ln\sigma_x(t_k)\rn)}\ge\eps.$$
Now,
$${\left|1-\bar{\sigma_j(t_k)}{1\over d_x}(\sigma(t_k),\check x)x_j
        \right|\over1-|\tilde\sigma_x(t_k)|^2}=\left|1-{\bar{\alpha_j(t_k)}
        \tilde\sigma_x(t_k)x_j\over1-|\tilde\sigma_x(t_k)|^2}\right|
        \le 1+{|\alpha_j(t_k)|\over1-|\tilde\sigma_x(t_k)|^2}.$$
Noticing that the function $t\mapsto t/(1+t)$ is strictly increasing we get
$$\left|{\sigma_j(t_k)-{1\over d_x}(\sigma(t_k),\check x)x_j\over
        1-\bar{\sigma_j(t_k)}{1\over d_x}(\sigma(t_k),\check x)x_j}\right|\ge
        {|\alpha_j(t_k)|/(1-|\tilde\sigma_x(t_k)|^2)\over1+
        |\alpha_j(t_k)|/(1-|\tilde\sigma_x(t_k)|^2)}\ge{\eps\over1+\eps}>0,$$
and so $\sigma$ is not special.\qedn

{\it Remark 1.1:} Instead of special $x$-curves, Jafari~[J] considers
$x$-curves {\sl tangent to the diagonal,} that is such that
$$\lim_{t\to 1^-}{\ln\sigma(t)-\sigma_x(t)\rn\over1-\ln\sigma_x(t)\rn}=0.$$
The previous proposition shows that if $x$ belongs to the \v Silov boundary
(the only case considered by Jafari) then a $x$-curve is special iff it is
tangent to the diagonal. But if~$x$ does not belong to the \v Silov boundary,
this is not true anymore. For instance, take~$x=(1,0)\in\de\Delta^2$, and
define $\sigma\colon[0,1)\to\Delta^2$ by
$$\sigma(t)=\bigl(\mezzo(1+t),\mezzo(1-t)\bigr).$$
Then it is easy to check that $\sigma$ is special but
$${\ln\sigma-\sigma_x\rn\over1-\ln\sigma_x\rn}\equiv 1,$$
and so $\sigma$ is not tangent to the diagonal in the sense of Jafari.
\smallbreak
We shall need another class of $x$-curves. We say that a $x$-curve
$\sigma$ is {\sl restricted} if $\tilde\sigma_x$ approaches~1 non-tangentially.
Notice that, again, being restricted imposes no conditions on the internal
components.

There is a more quantitative way of saying that a $x$-curve is restricted. The
{\sl Stolz region} $H(1,M)\subset\Delta$ of {\sl vertex}~$\tau\in\de\Delta$ and
{\sl amplitude}~$M>1$ is given by
$$H(\tau,M)=\left\{\zeta\in\Delta\biggm|{|\tau-\zeta|\over1-|\zeta|}<M\right
\}.$$
Geometrically, $H(\tau,M)$ is a sort of angle with vertex at~$\tau$, symmetric
with respect to the real axis, and amplitude going to~$\pi$ as~$M\to+\infty$;
therefore a 1-curve in~$\Delta$ approaches~1 non-tangentially iff it eventually
belongs to a Stolz region $H(1,M)$ for some~$M>1$. Thus we shall say that a
$x$-curve~$\sigma$ in~$\Delta^n$ is {\sl $M$-restricted} if~$\tilde\sigma_x$
eventually belongs to~$H(1,M)$.

Closely related to Stolz regions are the horocycles. The {\sl
horocycle}~$E(\tau,R)\subset\Delta$ of {\sl center}~$\tau\in\de\Delta$ and
{\sl radius}~$R>0$ is the set
$$E(\tau,R)=\left\{\zeta\in\Delta\biggm|{|\tau-\zeta|^2\over1-|\zeta|^2}<R
        \right\}.$$
Geometrically, $E(\tau,R)$ is an euclidean disk of euclidean radius~$R/(1+R)$
internally tangent to~$\de\Delta$ in~$\tau$.

In the following we shall need a version of Stolz regions and
horocycles in~$\Delta^n$. Take~$x\in\de\Delta^n$, $R>0$ and~$M>1$. Then the
({\sl
small\/}) {\sl horosphere}~$E(x,R)$ of {\sl center}~$x$ 
and {\sl radius}~$R$ and
the ({\sl small\/}) {\sl Kor\'anyi region}~$H(x,M)$ of {\sl
  vertex}~$x$ 
and {\sl
amplitude}~$M$ are defined by (see~[A1,~3]):
$$\displaylines{E(x,R)=\bigl\{z\in\Delta^n\bigm|\limsup\limits_{w\to
   x}\,[k_{\Delta^n}(z,w)-k_{\Delta^n}(0,w)]<\mezzo\log R\bigr\},\cr
 H(x,M)=\bigl\{z\in\Delta^n\bigm|\limsup\limits_{w\to x}\,[k_{\Delta^n}(z,w)-
        k_{\Delta^n}(0,w)]+k_{\Delta^n}(0,z)<\log M\bigr\}.\cr}$$
It is easy to check that if $x=\tau\in\de\Delta$ then we recover the horocycles
and Stolz regions previously defined.
\smallbreak
{\it Remark 1.2:} Replacing the $\limsup$ by $\liminf$ one gets the definitions
of big horosphere~$F(x,R)$ and big Kor\'anyi region~$K(x,M)$, but we shall not
need them in this paper.
\smallbreak

To understand the shape of horospheres and Kor\'anyi regions in the polydisk we
need to compute the $\limsup$ used to define them:

\newthm Proposition \Lim: Given $x\in\de\Delta^n$ and $z\in\Delta^n$ one has
$$\eqalign{\limsup_{w\to x}[k_{\Delta^n}(z,w)-k_{\Delta^n}(0,w)]&=\mezzo\log
        \max_{|x_j|=1}\left\{{|x_j-z_j|^2\over1-|z_j|^2}\right\}\cr
        &=\lim_{s\to1^-}[k_{\Delta^n}\bigl(z,\phe_x(s)\bigr)-\omega(0,s)].\cr}
        \neweq\eqlim$$

\pf First of all,
$$k_{\Delta^n}(z,w)-k_{\Delta^n}(0,w)=\log\left({1+\ln\gamma_z(w)\rn\over1+\ln
        w\rn}\right)
+\mezzo\log\left({1-\ln w\rn^2\over1-\ln\gamma_z(w)\rn^2}\right).$$
Since $\ln\gamma_z(x)\rn=\ln x\rn=1$, we need to study the behavior of the
second addend only. Now
$$\displaylines{1-\ln w\rn^2=\min_h\{1-|w_h|^2\};\cr
        1-\ln\gamma_z(w)\rn^2
=\min_j\left\{{1-|z_j|^2\over|1-\bar{z_j}w_j|^2}(1-
        |w_j|^2)\right\}.\cr}$$
Hence
$${1-\ln w\rn^2\over1-\ln\gamma_z(w)\rn^2}
=\max_j\left\{{|1-\bar{z_j}w_j|^2\over
        1-|z_j|^2}\,\min_h\left\{{1-|w_h|^2\over1-|w_j|^2}\right\}\right\}.$$
If $|x_j|<1$ we have
$$\min_h\left\{{1-|w_h|^2\over1-|w_j|^2}\right\}\to0$$
as $w\to x$; so we ought to consider only $j$'s such that $x_j$ is a \v Silov
component of~$x$.

We clearly have
$$\min_h\left\{{1-|w_h|^2\over1-|w_j|^2}\right\}\le1;$$
therefore
$$\limsup_{w\to x}{1-\ln w\rn^2\over1-\ln\gamma_z(w)\rn^2}\le\max_{|x_j|=1}
        \left\{{|1-\bar{z_j}w_j|^2\over1-|z_j|^2}\right\}.\neweq\eqkh$$
To prove the opposite inequality, set $w^\nu=(1-1/\nu)^{1/2}x$; we have
$$1-|(w^\nu)_h|^2=(1-|x_h|^2)+|x_h|^2/\nu.$$
Hence if $|x_j|=1$ we get
$$\lim_{\nu\to\infty}\min_h\left\{{1-|(w^\nu)_h|^2\over1-|(w^\nu)_j|^2}\right\}
        =1,$$
and we have proved the first equality in~\eqlim.

For the second equality, notice that
$$k_{\Delta^n}\bigl(z,\phe_x(s)\bigr)-\omega(0,s)=\log\left({1+\ln\gamma_z
        \bigl(\phe_x(s)\bigr)\rn\over1+s}\right)+\mezzo\log\left({1-s^2\over
        1-\ln\gamma_z\bigl(\phe_x(s)\bigr)\rn^2}\right).$$
The first addend on the right goes to~0 as~$s\to 1^-$. Next we have
$$1-\ln\gamma_z\bigl(\phe_x(s)\bigr)\rn^2=\min_j\left\{{1-|z_j|^2\over
        |1-s\bar{z_j}x_j|^2}(1-s^2|x_j|^2)\right\};$$
therefore
$${1-s^2\over1-\ln\gamma_z\bigl(\phe_x(s)\bigr)\rn^2}=\max_j\left\{{|1-s\bar
{z_j}
        x_j|^2\over1-|z_j|^2}\,{1-s^2\over1-s^2|x_j|^2}\right\}.$$
Since $(1-s^2)/(1-s^2|x_j|^2)\equiv1$ if $|x_j|=1$ and it goes to~0 if
$|x_j|<1$, we are done.\qedn

Therefore for every $x\in\de\Delta^n$, $R>0$ and $M>1$ we have
$$E(x,R)=\left\{z\in\Delta^n\biggm|\max_{|x_j|=1}\left\{{|x_j-z_j|^2\over
        1-|z_j|^2}\right\}<R\right\},$$
and
$$H(x,M)=\left\{z\in\Delta^n\biggm|{1+\ln z\rn\over1-\ln z\rn}\,\max_{|x_j|=1}
        \left\{{|x_j-z_j|^2\over1-|z_j|^2}\right\}<M^2\right\}.\neweq\eqHxR$$
The shape of $E(x,R)$ is easily described: it is a product $E_1\times\cdots
\times E_n$, where $E_j=E(x_j,R)$ is a horocycle in~$\Delta$ if $|x_j|=1$, and
$E_j=\Delta$ otherwise.

The shape of $H(x,M)$ is more complicated, but anyway we can prove the
following:

\newthm Proposition \shape: Let $x\in\de\Delta^n$ and $M>1$. Then
$$\bigcup_{t\ge0} B_{\Delta^n}\bigl(\phe_x(t),\mezzo\log M)\subseteq H(x,M)
        \subseteq H_1\times\cdots\times H_n,\neweq\eqinc$$
where $B_{\Delta^n}(z,r)$ is the open ball of center~$z$ and radius~$r$ with
respect to the Kobayashi distance of~$\Delta^n$,
$H_j=H(x_j,M)$ is a Stolz region in~$\Delta$ if $|x_j|=1$ and~$H_j=\Delta$
otherwise.

\pf For every $t\ge0$ we have
$$\lim_{s\to1^-}[k_{\Delta^n}\bigl(z,\phe_x(s)\bigr)-\omega(0,s)]+
        k_{\Delta^n}(0,z)\le 2k_{\Delta^n}\bigl(z,\phe_x(t)\bigr),$$
because $k_{\Delta^n}\bigl(\phe_x(t),\phe_x(s)\bigr)-\omega(0,s)+k_{\Delta^n}
\bigl(0,\phe_x(t)\bigr)=0$ as soon as~$s\ge t$. Proposition~\Lim\ then implies
the first inclusion in~\eqinc.

For the second inclusion, take $z\in H(x,M)$ and suppose that $|x_h|=1$. Then
$${|x_h-z_h|^2\over(1-|z_h|)^2}={1+|z_h|\over1-|z_h|}\,{|x_h-z_h|^2\over1-
        |z_h|^2}\le{1+\ln z\rn\over1-\ln z\rn}\,\max_{|x_j|=1}\left\{
        {|x_j-z_j|^2\over1-|z_j|^2}\right\}<M^2,$$
and so $z_j\in H(x_j,M)$.\qedn

{\it Remark 1.3:} Given $x\in\de\Delta^n$, let $z=t\check x+v$, where
$t\in[0,1)$, the vector~$v$ has no \v Silov components with respect to~$x$ and
$\ln v\rn\le t$. Then it is easy to check that $z\in H(x,M)$ for all~$M>1$; in
particular, in the second inclusion of~\eqinc\ we are forced to take the whole
unit disk as factor for the internal components.
\smallbreak
{\it Remark 1.4:} The second inclusion in~\eqinc\ implies that if $z\to x$
inside some $H(x,M)$ and $x_j$ is a \v Silov component of~$x$, 
then $z_j\to x_j$
non-tangentially; on the other hand, the first inclusion implies that $z_j\to
x_j$ unrestricted if $x_j$ is an internal component.
\smallbreak
{\it Remark 1.5:} In~[Kr] Krantz suggested to define approach regions
in general domains as unions like the one in the left-hand member of \eqinc,
replacing the geodesic line~$\phe_x(t)$ by the euclidean normal line at~$x$.
Proposition~\shape\ shows that our Kor\'anyi regions are comparable with
Krantz's approach regions.
\smallbreak
Another consequence of Proposition~\Lim\ is the following:

\newthm Corollary \add: Take $x\in\de\Delta^n$, $R>0$ and $M>1$. Then
$$\phe_x\bigl(E(1,R)\bigr)=E(x,R)\cap\phe_x(\Delta)=p_x\bigl(E(x,R)\bigr)$$
and
$$\phe_x\bigl(H(1,M)\bigr)=H(x,M)\cap\phe_x(\Delta)=p_x\bigl(
        H(x,M)\bigr).$$

\pf Since $\phe_x$ is a complex geodesic, we have
$k_{\Delta^n}\bigl(\phe_x(\zeta),\phe_x(s)\bigr)=\omega(\zeta,s)$ for
every~$\zeta\in\Delta$ and~$s\in[0,1)$; the assertions then follows from
Proposition~\Lim.\qedn

If $n=1$ we have seen that there is a strong relationship between
non-tangential curves and Stolz regions. A similar fact holds for~$n>1$ too:

\newthm Proposition \rcKr: 
Take $x\in\de\Delta^n$, and let $\sigma\colon[0,1)\to
\Delta^n$ be a $x$-curve. Fix~$M>1$. Then:
{\smallskip
\itm{(i)}$\sigma$ is $M$-restricted iff $\sigma_x(t)\in H(x,M)$ eventually;
\itm{(ii)}if $\sigma(t)\in H(x,M)$ eventually then $\sigma$ is $M$-restricted;
\itm{(iii)}if $\sigma$ is special and $M$-restricted then for any $M_1>M$ we
have $\sigma(t)\in H(x,M_1)$ eventually.}

\pf (i) By definition, $\sigma$ is $M$-restricted iff $\tilde\sigma_x(t)\in
H(1,M)$ eventually iff
$$\lim_{s\to 1^-}[\omega(\tilde\sigma_x(t),s)-\omega(0,s)]+\omega\bigl(0,
        \tilde\sigma_x(t)\bigr)<\log M$$
eventually. Being $\phe_x$ a complex geodesic, and since
$\sigma_x=\phe_x\circ\tilde\sigma_x$, this happens iff eventually
$$\lim_{s\to1^-}[k_{\Delta^n}\bigl(\sigma_x(t),\phe_x(s)\bigr)-\omega(0,s)]+
        k_{\Delta^n}\bigl(0,\sigma_x(t)\bigr)<\log M,$$
that is iff $\sigma_x(t)\in H(x,M)$ eventually, by Proposition~\Lim.
\smallbreak
(ii) For any $z\in\Delta^n$ we have
$$\eqalign{\lim_{s\to1^-}[k_{\Delta^n}\bigl(p_x(z),\phe_x(s)\bigr)-\omega(0,
s)]&+
        k_{\Delta^n}\bigl(0,p_x(z)\bigr)\cr
        &\le\lim_{s\to1^-}[k_{\Delta^n}\bigl(z,\phe_x(s)\bigr)-\omega(0,s)]
        +k_{\Delta^n}(0,z),\cr}$$
because $p_x\circ\phe_x=\phe_x$, and the assertion follows, again by
Proposition~\Lim.
\smallbreak
(iii) We have
$$\eqalign{k_{\Delta^n}
\bigl(\sigma(t),\phe_x(s)\bigr)&-\omega(0,s)+k_{\Delta^n}
    \bigl(0,\sigma(t)\bigr)\cr
  &\le2\kappa_{\Delta^n}\bigl(\sigma(t),\sigma_x(t)\bigr)+k_{\Delta^n}\bigl(
  \sigma_x(t),\phe_x(s)\bigr)-\omega(0,s)+\kappa_{\Delta^n}\bigl(0,\sigma_x(t)
        \bigr),\cr}$$
and the assertion follows from (i) and Proposition~\Lim.\qedn

{\it Remark 1.6:} There are restricted $x$-curves that does not
eventually belong to any Kor\'anyi region. For instance, take~$x=(1,1)\in\de
\Delta^2$, and let $\sigma\colon[0,1)\to\Delta^2$ be given by
$$\sigma(t)=\bigl(t+i(1-t)^{1/2},t-i(1-t)^{1/2}\bigr).$$
Clearly $\tilde\sigma_x(t)=t$, and so $\sigma$ is $M$-restricted for any~$M>1$.
On the other hand,
$$1-\ln\sigma(t)\rn=O(1-t)=1-|\sigma_j(t)|^2,\qquad |1-\sigma_j(t)|^2=O(1-t),$$
and so
$${1+\ln\sigma(t)\rn\over1-\ln\sigma(t)\rn}\,\max_j\left\{{|1-\sigma_j(t)|^2
        \over1-|\sigma_j(t)|^2}\right\}=O\bigl((1-t)^{-1}\bigr)\to+\infty,$$
as claimed. Notice that
$${\ln\sigma(t)-\sigma_x(t)\rn\over1-\ln\sigma_x(t)\rn}={1\over(1-t)^{1/2}}
        \to+\infty,$$
and $\sigma$ is not special.
\smallbreak
We end this section with a result we shall use to build examples later on:

\newthm Lemma \fun: Let $\alpha_1$,~$\alpha_2\in\Hol(\Delta^2,\C)$ be so that
$\Re\alpha_j(z)>|\Im\alpha_j(z)|$ for all $z\in\Delta^2$ and~$j=1$,~$2$.
Then the
function
$$f(z)={\alpha_1(z)-\alpha_2(z)\over\alpha_1(z)+\alpha_2(z)}$$
is a holomorphic function from~$\Delta^2$ into~$\Delta$.

\pf The hypothesis ensures that $\Re\bigl(\alpha_1(z)\bar{\alpha_2(z)}\bigr)>0$
for all~$z\in\Delta^2$. Therefore
$$\forevery{z\in\Delta^2}|\alpha_1(z)+\alpha_2(z)|^2>|\alpha_1(z)-\alpha_2(z)|
        ^2,$$
and we are done.\qedn

\smallsect 2. The Lindel\"of principle

In the previous section we have defined, via the canonical projection device,
special and restricted $x$-curves. We say that a map
$f\colon\Delta^n\to\C^m$ has {\sl restricted
$K$-limit}~$L\in\C^m$ at~$x\in\de\Delta^n$ if
$f\bigl(\sigma(t)\bigr)\to L$ as~$t\to1^-$ for any special restricted
$x$-curve~$\sigma$; we shall write
$$\rKlim_{z\to x}f(z)=L.$$
We say that $f$ has {\sl $K$-limit}~$L\in\C^m$
at~$x$ if~$f(z)\to L$ as~$z\to x$ inside any Kor\'anyi region~$H(x,M)$.
\smallbreak
{\it Remark 2.1:} By Proposition~\rcKr.(iii), if $f$ has $K$-limit~$L$ at~$x$
then it has restricted $K$-limit~$L$ at~$x$ too. 
The converse is false, even for
bounded holomorphic functions: let $f\colon\Delta^2\to\C$ be given by
$$f(z_1,z_2)={(1-z_1)^{1/2}-(1-z_2)^{1/2}\over(1-z_1)^{1/2}+(1-z_2)^{1/2}}.$$
By Lemma~\fun, the image of~$f$ is contained in~$\Delta$. Set $x=(1,1)\in\de
\Delta^2$, and take a special restricted $x$-curve~$\sigma$. Write
$\sigma=\sigma_x+\alpha$; then
$$f(\sigma)
={(1-\tilde\sigma_x-\alpha_1)^{1/2}-(1-\tilde\sigma_x-\alpha_2)^{1/2}
   \over(1-\tilde\sigma_x-\alpha_1)^{1/2}+(1-\tilde\sigma_x-\alpha_2)^{1/2}}=
   {\left(1-{\alpha_1\over1-\tilde\sigma_x}\right)^{1/2}-\left(1-{\alpha_2\over
  1-\tilde\sigma_x}\right)^{1/2}\over\left(1-{\alpha_1\over1-\tilde\sigma_x}
        \right)^{1/2}+\left(1-{\alpha_2\over1-\tilde\sigma_x}\right)^{1/2}}.$$
Being $\sigma$ special we have
$$\left|{\alpha_j\over1-\tilde\sigma_x}\right|\le{\ln\alpha\rn\over
        1-|\tilde\sigma_x|}={\ln\sigma-\sigma_x\rn\over1-\ln\sigma_x\rn}\to0,$$
and so $f$ has restricted $K$-limit 0 at~$x$.

On the other hand, for $\lambda\in(0,1)$ let $\sigma^\lambda\colon[0,1)\to
\Delta^2$ be the $x$-curve
$$\sigma^\lambda(t)=\bigl(t,t+\lambda(1-t)\bigr).$$
Then
$${1+\ln\sigma^\lambda\rn\over1-\ln\sigma^\lambda\rn}\max_j\left\{
        {|1-\sigma^\lambda_j|^2\over1-|\sigma^\lambda_j|^2}\right\}
={(1+\lambda)+
        (1-\lambda)t\over(1+t)(1-\lambda)}\le{2\over1-\lambda},$$
that is $\sigma^\lambda(t)\in H\bigl(x,2/(1-\lambda)\bigr)$. But $f\bigl(
\sigma^\lambda(t)\bigr)\equiv\bigl(1-(1-\lambda)^{1/2}\bigr)\big/\bigl(1+(1-
\lambda)^{1/2}\bigr)$, and so~$f$ has no $K$-limit at~$x$.
\smallbreak
In the next section we shall use yet again another kind of limit, stronger than
$K$-limit. Take~$x\in\de\Delta^n$; we shall say that a $x$-curve~$\sigma$ is
{\sl peculiar} if~$\sigma(t)\in E(x,R)$ eventually for
all~$R>0$. 
Recalling the
shape of horospheres, this means that
$$\lim_{t\to1^-}\max_{|x_j|=1}
\left\{{|x_j-\sigma_j(t)|^2\over1-|\sigma_j(t)|^2}
        \right\}=0;\neweq\eqElim$$
as usual, being peculiar imposes no restrictions on the internal components.

We say that a function $f\colon\Delta^n\to\C^m$ admits
{\sl restricted $E$-limit}~$L\in\C^m$ at~$x\in\de\Delta^n$ if $f\bigl(\sigma(t)
\bigr)\to L$ as~$t\to1^-$ for any peculiar $x$-curve~$\sigma$.
\smallbreak
{\it Remark 2.2:} One could obviously say that a
function~$f\colon\Delta^n\to\C^m$ has {\sl $E$-limit}~$L$ at~$x$ if~$f(z)\to L$
as~$z\to x$ inside any horosphere~$E(x,R)$, but we shall not use this
definition.
\smallbreak
{\it Remark 2.3:} It is easy to check that for
every~$M>1$ and~$R>0$ one has
$$H(x,M)\setminus B_{\Delta^n}(0,r)\subset E(x,R),$$
where $r=\mezzo\log(M^2/R)$; it follows that if~$f$ has restricted
$E$-limit~$L$ at~$x$ then it has \hbox{$K$-limit}~$L$ there. To prove that the
converse does not hold, not even for bounded holomorphic functions, we need a
couple of preliminary observations.

First of all, take $z\in H\bigl((1,1),M)$. Then
$${1\over2M^2}\le{1-|z_2|\over1-|z_1|}\le2M^2.\neweq\eqcku$$
Indeed we have
$$\eqalign{M^2\ge{1+\ln z\rn\over1-\ln z\rn}\max_j\left\{{|1-z_j|^2\over
   1-|z_j|^2}\right\}&\ge{1+|z_1|\over1-|z_1|}\,{|1-z_2|^2\over1-|z_2|^2}\cr
        &\ge{1+|z_1|\over1-|z_1|}\,{(1-|z_2|)^2\over1-|z_2|^2}\ge{1\over2}\,
        {1-|z_2|\over1-|z_1|},\cr}\neweq\eqckm$$
and the right inequality in~\eqcku\ follows; the other one is obtained
reversing the roles of~$z_1$ and~$z_2$. As a consequence,
$${1\over2M^2}\le\left|{1-z_2\over1-z_1}\right|\le2M^2.\neweq\eqckd$$
Indeed \eqcku\ and \eqckm\ imply
$$\left|{1-z_2\over1-z_1}\right|^2\le{|1-z_2|^2\over(1-|z_1|)^2}\le
M^2{1-|z_2|^2
        \over1-|z_1|^2}\le4M^4$$
for all $z\in H\bigl((1,1),M\bigr)$.

Now fix $\alpha<1$ and define $f\in\Hol(\Delta^2,\C)$ by
$$f(z_1,z_2)={(1-z_1)^{\alpha/2}-(1-z_2)^{1/2}\over(1-z_1)^{\alpha/2}+(1-z_2)
        ^{1/2}};$$
Lemma~\fun\ ensures us that $f(\Delta^2)\subset\Delta$. 
Since $\alpha<1$ we have
$$f(z_1,z_2)={1-\left({1-z_2\over1-z_1}\right)^{\alpha/2}(1-z_2)^{(1-\alpha)/2}
 \over1+\left({1-z_2\over1-z_1}\right)^{\alpha/2}(1-z_2)^{(1-\alpha)/2}},$$
and so, by \eqckd, $f$ has $K$-limit~$1$ at~$(1,1)$. On the other hand, for
$\lambda\in(0,1)$ let $\sigma^\lambda$ be the $(1,1)$-curve given by
$$\sigma^\lambda(t)=\bigl(t,t-\lambda(1-t)^\alpha\bigr).$$
It is easy to check that $\sigma^\lambda$ is peculiar, but
$$f\bigl(\sigma^\lambda(t)\bigr)={(1-t)^{\alpha/2}-[(1-t)+\lambda(1-t)^\alpha]
   ^{1/2}\over (1-t)^{\alpha/2}+[(1-t)+\lambda(1-t)^\alpha]^{1/2}}
  ={1-[\lambda+(1-t)^{1-\alpha}]^{1/2}\over1+[\lambda+(1-t)^{1-\alpha}]^{1/2}}
        \to{1-\lambda^{1/2}\over1+\lambda^{1/2}},$$
and so $f$ has no restricted $E$-limit at~$(1,1)$.
\smallbreak
{\it Remark 2.4:} The classical Lindel\"of principle implies that
examples like the previous one cannot exist in~$\Delta$. Indeed, if
$f\in\Hol(\Delta,\Delta)$ has limit~$L$ along any given 1-curve, then~$L$ is
the non-tangential limit of~$f$ at~$1$; therefore if $f$ restricted to
any other 1-curve admits a limit, that limit should be~$L$.
\smallbreak
Using these definitions it is very easy to prove a Lindel\"of principle
(Theorem~\zLA):

\newthm Theorem \Lin: Let $f\colon\Delta^n\to\C^m$ be a bounded holomorphic
function. Given $x\in\de\Delta^n$, assume there is a special
$x$-curve~$\sigma^o$ such that
$$\lim_{t\to1^-}f\bigl(\sigma^o(t)\bigr)=L\in\C^m.$$
Then $f$ has restricted $K$-limit~$L$ at~$x$.

\pf (see~[A3]). Clearly we can assume $m=1$ and
$f(\Delta^n)\subset\subset\Delta$. Let $\sigma$ be any special $x$-curve. We
have
$$\omega\left(f\bigl(\sigma(t)\bigr),f\bigl(\sigma_x(t)\bigr)\right)\le
        k_{\Delta^n}\bigl(\sigma(t),\sigma_x(t)\bigr)\to0;$$
therefore the limit of $f\bigl(\sigma(t)\bigr)$ as $t\to1^-$ exists iff the
limit of~$f\bigl(\sigma_x(t)\bigr)$ as~$t\to1^-$ does, and the two are equal.
In particular, $f\bigl(\sigma^o_x(t)\bigr)\to L$ as~$t\to1^-$. Hence, by the
classical Lindel\"of principle (see~[R]), $f\bigl(\sigma_x(t)\bigr)\to L$ for
any restricted $x$-curve~$\sigma$ and, by the previous observation, $f\bigl(
\sigma(t)\bigr)\to L$ for any special restricted $x$-curve~$\sigma$.\qedn

{\it Remark 2.5:} The bounded holomorphic functions~$f\in\Hol(\Delta^2,\Delta)$
described in Remarks~2.1 and~2.3 show that if $f\bigl(\sigma(t)\bigr)\to L$
as~$t\to1^-$ but~$\sigma$ is not special then it is not necessarily true that
$f$ has restricted $K$-limit~$L$.
\smallbreak
It turns out that to prove a Julia-Wolff-Carath\'eodory theorem a stronger
result is needed. Given~$x\in\de\Delta^n$, if $f\colon\Delta^n\to\C^m$ is such
that for every~$M>1$ there is a constant~$c_M>0$ such that $\|f(z)\|\le
c_M$ for all~$z\in H(x,M)$, we shall say that $f$ is {\sl $K$-bounded} at~$x$.
Then:

\newthm Theorem \LinK: Given $x\in\de\Delta^n$, let $f\colon\Delta^n\to\C^m$ be
a holomorphic map $K$-bounded at~$x$. Assume there is a restricted special
$x$-curve~$\sigma^o$ such that
$$\lim_{t\to1^-}f\bigl(\sigma^o(t)\bigr)=L\in\C^m.$$
Then $f$ has restricted $K$-limit $L$ at~$x$.

\pf We can of course assume $m=1$. First of all we claim that if $\sigma$ is an
$M$-restricted special $x$-curve, then for all~$M_1>M$ we have
$$\lim_{t\to1^-}k_{H(x,M_1)}\bigl(\sigma(t),\sigma_x(t)\bigr)=0.\neweq\eqlimu$$
For any $t\in[0,1)$ let us consider the map $\psi_t\colon\C\to\C^n$ given by
$$\psi_t(\zeta)=\sigma_x(t)+\zeta\bigl(\sigma(t)-\sigma_x(t)\bigr).$$
Clearly, $\psi_t(0)=\sigma_x(t)$ and $\psi_t(1)=\sigma(t)$. Assume for the
moment that we have proved that for every~$R>1$ there
exists~$t_0=t_0(R)\in[0,1)$ such that
$$\forevery{t_0<t<1}\psi_t(\Delta_R)\subset H(x,M_1),\neweq\eqinc$$
where $\Delta_R=\{\zeta\in\C\mid|\zeta|<R\}$. Set
$$R(t)=\sup\{r>0\mid\psi_t(\Delta_r)\subset H(x,M_1)\}.$$
Inclusion \eqinc\ implies that $R(t)\to+\infty$ as~$t\to1^-$; since, by
definition,
$$k_{H(x,M_1)}\bigl(\sigma(t),\sigma_x(t)\bigr)\le\inf\left\{{1\over R}\biggm|
 \exists\phe\in\Hol\bigl(\Delta_R,H(x,M_1)\bigr):\phe(0)=\sigma_x(t),\phe(1)
        =\sigma(t)\right\},$$
equation \eqlimu\ follows from \eqinc.

Now we prove \eqinc. If we write $\sigma=\tilde\sigma_xx
+\alpha$, clearly we have~$\psi_t(\zeta)=\tilde\sigma_x(t)x+\zeta\alpha(t)$,
with~$\alpha(t)\to0$ as~$t\to1^-$. In particular, for every~$R>1$ there
exists~$t_1=t_1(R)\in[0,1)$ such that
$$\forevery{\zeta\in\Delta_R\quad\forall t_1<t<1}\quad
        \qquad\qquad\ln\psi_t(\zeta)\rn=\max_{|x_h|=1}
        \bigl\{|\tilde\sigma_x(t)x_h+\zeta\alpha_h(t)|\bigr\}.\neweq\eqpass$$
Assume, by contradiction, that \eqinc\ is false. Then there exist $M_1>M$
and~$R_0>1$ such that for any~$t_0\in[0,1)$ there are $t'=t'(t_0)\in(t_0,1)$
and~$\zeta_0=\zeta_0(t_0)\in\Delta_{R_0}$ such that $\psi_{t'}(\zeta_0)\notin
H(x,M_1)$. Since $\sigma_x(t')\in H(x,M_1)$ eventually (Proposition~\rcKr), we
can choose~$t_1=t_1(R_0)\in(0,1)$ such that $\sigma_x(t')\in H(x,M_1)$ for
all~$t_0>t_1$; being~$H(x,M_1)$ open we can also assume that
$\psi_{t'}(\zeta_0)\in\de H(x,M_1)$ but $\psi_{t'}(\zeta)\in H(x,M_1)$ for
all~$\zeta\in\Delta_{|\zeta_0|}$.

Recalling \eqHxR\ and \eqpass\ we get
$$\max_{|x_h|=1}\left\{{1+|\tilde\sigma_xx_h
+\zeta_0\alpha_h|\over1-|\tilde\sigma_xx_h
+\zeta_0\alpha_h|}\right\}\,\max_{|x_j|=1}\left\{{|1-\tilde\sigma_x-
    \zeta_0\alpha_j\bar{x_j}|^2\over1-|\tilde\sigma_xx_j+\zeta_0\alpha_j|^2}
\right\}=M^2_1,$$ where everything is evaluated at~$t=t'$. So
there are, possibly different, \v Silov components~$x_{h_0}$
and~$x_{j_0}$ such that
$$M^2_1=\left({|1-\tilde\sigma_x-\zeta_0\alpha_{j_0}\bar{x_{j_0}}|\over
        1-|\tilde\sigma_x+\zeta_0\alpha_{j_0}\bar{x_{j_0}}|}\right)^2\left[
        {1+|\tilde\sigma_xx_{h_0}+\zeta_0\alpha_{h_0}|\over1+
        |\tilde\sigma_xx_{j_0}+\zeta_0\alpha_{j_0}|}\bigg/{1-
        |\tilde\sigma_xx_{h_0}+\zeta_0\alpha_{h_0}|\over1-
        |\tilde\sigma_xx_{j_0}+\zeta_0\alpha_{j_0}|}\right].\neweq\eqdf$$
Now, if $\zeta\in\Delta_{R_0}$ we have
$$\eqalign{{1+|\zeta|{\displaystyle{|\alpha_{j_0}|\over1-
|\tilde\sigma_x|}}\over
        1-|\zeta|{\displaystyle{|\alpha_{h_0}|\over1-|\tilde\sigma_x|}}}=
        {1-|\tilde\sigma_x|+|\zeta||\alpha_{j_0}|\over1-|\tilde\sigma_x|-
 |\zeta||\alpha_{h_0}|}&\ge{1-|\tilde\sigma_xx_{j_0}+\zeta\alpha_{j_0}|\over
        1-|\tilde\sigma_xx_{h_0}+\zeta\alpha_{h_0}|}\cr
        &\ge{1-|\tilde\sigma_x|-|\zeta||\alpha_{j_0}|\over1-|\tilde\sigma_x|+
        |\zeta||\alpha_{h_0}|}={1-|\zeta|{\displaystyle{|\alpha_{j_0}|\over
        1-|\tilde\sigma_x|}}\over1+|\zeta|{\displaystyle{|\alpha_{h_0}|\over
        1-|\tilde\sigma_x|}}}.\cr}\neweq\eqcp$$
Being $\sigma$ special, Proposition~\Jafari\ yields
$$\max_{|x_j|=1}{|\alpha_j|\over1-|\tilde\sigma_x|}\to0;\neweq\eqJafu$$
therefore \eqcp\ implies that
$${1-|\tilde\sigma_xx_{j_0}+\zeta\alpha_{j_0}|\over
        1-|\tilde\sigma_xx_{h_0}+\zeta\alpha_{h_0}|}\to1\neweq\eqlind$$
uniformly for $\zeta\in\Delta_{R_0}$. Now fix $\eps>0$ so that $M'_1=M_1/(1+
\eps)>M$; by~\eqlind\ we can choose a sequence~$t^k_0\to1^-$ so that~\eqdf\
holds for all~$t^k_0$ with the same~$j_0$ and~$h_0$, and moreover
$$\left[{1+|\tilde\sigma_xx_{h_0}+\zeta\alpha_{h_0}|\over1+
        |\tilde\sigma_xx_{j_0}+\zeta\alpha_{j_0}|}\bigg/{1-
        |\tilde\sigma_xx_{h_0}+\zeta\alpha_{h_0}|\over1-
        |\tilde\sigma_xx_{j_0}+\zeta\alpha_{j_0}|}\right]\le(1+\eps)^2$$
for all $\zeta\in\Delta_{R_0}$ and all~$t^k_0$ (where everything is evaluated
at~$t'(t^k_0)$, as usual). Recalling~\eqdf\ we then get
$${|1-\tilde\sigma_x-\zeta_0\alpha_{j_0}\bar{x_{j_0}}|\over
        1-|\tilde\sigma_x+\zeta_0\alpha_{j_0}\bar{x_{j_0}}|}\ge{M_1\over1+\eps}
        =M'_1>M,$$
again for all~$t^k_0$. Writing $v_k=\zeta_0(t^k_0)\alpha_{j_0}\bigl(t'(t^k_0)
\bigr)\bar{x_{j_0}}$ we obtain
$$|1-\tilde\sigma_x|+|v_k|\ge|1-\tilde\sigma_x-v_k|\ge M'_1(1-|\tilde\sigma_x
        +v_k|)\ge M'_1(1-|\tilde\sigma_x|-|v_k|);$$
therefore, being $\sigma$ $M$-restricted, for $k$ large enough we have
$$M+{|v_k|\over1-|\tilde\sigma_x|}
\ge{|1-\tilde\sigma_x|\over1-|\tilde\sigma_x|}
        +{|v_k|\over1-|\tilde\sigma_x|}
\ge M'_1\left(1-{|v_k|\over1-|\tilde\sigma_x|}
        \right),$$
whence
$${|v_k|\over1-|\tilde\sigma_x|}\ge{M'_1-M\over1+M'_1},$$
and so
$$R_0\ge{M'_1-M\over1+M'_1}\,{1-\bigl|\tilde\sigma_x\bigl(t'(t^k_0)\bigr)\bigr|
        \over\bigl|\alpha_{j_0}\bigl(t'(t^k_0)\bigr)\bigr|}.$$
Letting $k\to+\infty$, that is letting $t'(t^k_0)\to1^-$, we finally get a
contradiction, because of~\eqJafu.

Summing up, we have proved that \eqinc\ holds, and therefore \eqlimu\ holds
too; we can now finish the proof of the theorem. Let $M>1$ so that $\sigma^o$
is~$M$-restricted, and fix~$M_1>M$; then \eqlimu\ holds. On~$H(x,M_1)$ the
function~$f$ is bounded by~$c$, say; therefore
$$k_{\Delta_c}\left(f\bigl(\sigma^o(t)\bigr),f\bigl(\sigma^o_x(t)\bigr)\right)
        \le k_{H(x,M_1)}\bigl(\sigma^o(t),\sigma^o_x(t)\bigr),$$
and so
$$\lim_{t\to1^-}f\bigl(\sigma^o_x(t)\bigr)=L. \neweq\eqlimt$$
Finally, let $\sigma$ be any restricted special $x$-curve. The classical
Lindel\"of principle applied to~$f\circ\phe_x$ together with \eqlimt\ implies
$$\lim_{t\to1^-}f\bigl(\sigma_x(t)\bigr)=L;$$
hence, arguing as before, we find that $f\bigl(\sigma(t)\bigr)\to L$ as~$t\to
1^-$, and we are done.\qedn

\smallsect 3. Julia's lemma

Now we are ready to deal with Julia's lemma in the polydisk, at least for
functions; the general case is discussed in the last section (but see
also Theorem~\zJA).

The standard Julia's lemma in the disk says that if $f\colon\Delta\to\Delta$ is
a bounded holomorphic function such that the rate of approach of~$f(\zeta)$
to~$\de\Delta$ is comparable to the rate of approach of~$\zeta$
to~$\sigma\in\de\Delta$, then $f$ sends horocycles centered at~$\sigma$ into
horocycles centered at some~$\tau\in\de\Delta$ --- and then $f$ has
non-tangential limit~$\tau$ at~$\sigma$.

If $f\in\Hol(\Delta^n,\Delta)$ and we want a version of Julia's lemma in the
polydisk, the natural thing to do is to compare the rate of approach of~$f(w)$
to~$\de\Delta^n$ with the rate of approach of~$w$ to~$x\in\de\Delta^n$, that is
to study
$${1-|f(w)|\over1-\ln w\rn}$$
when $w\to x$. Now, it is easy to check that
$$\mezzo\log\liminf_{w\to x}{1-|f(w)|\over1-\ln w\rn}=\liminf_{w\to x}\bigl[
        k_{\Delta^n}(0,w)-\omega\bigl(0,f(w)\bigr)\bigr];\neweq\eqliminf$$
since we have defined horospheres and the like in terms of the Kobayashi
distance, the natural statement for a Julia lemma is:

\newthm Theorem \Julia: Let $f\colon\Delta^n\to\Delta$ be a bounded holomorphic
function, and let $x\in\de\Delta^n$ be such that
$$\liminf_{w\to x}\bigl[k_{\Delta^n}(0,w)-\omega\bigl(0,f(w)\bigr)\bigr]\le
        \mezzo\log\alpha<+\infty.\neweq\eqaJ$$
Then there exists $\tau\in\de\Delta$ such that
$$\forevery{R>0}f\bigl(E(x,R)\bigr)\subseteq E(\tau,\alpha R).\neweq\eqbJ$$
Furthermore, $f$ admits restricted $E$-limit~$\tau$ at~$x$.

\pf First of all choose a sequence $\{w_\nu\}\subset\Delta^n$ converging to~$x$
such that
$$\lim_{\nu\to\infty}\bigl[k_{\Delta^n}(0,w_\nu)-\omega\bigl(0,f(w_\nu)\bigr)
        \bigr]=\liminf_{w\to x}\bigl[k_{\Delta^n}(0,w)-\omega\bigl(0,f(w)\bigr)
        \bigr].$$
Up to a subsequence, we can assume that $f(w_\nu)\to\tau\in\bar\Delta$. Since
$\Delta^n$ is complete hyperbolic, we have~$k_{\Delta^n}(0,w_\nu)\to+\infty$;
therefore $\omega\bigl(0,f(w_\nu)\bigr)\to+\infty$ as well,
and~$\tau\in\de\Delta$.

Now take $z\in E(x,R)$; then
$$\eqalign{\lim_{\zeta\to\tau}[\omega(f(z),\zeta)&-\omega(0,\zeta)]
=\lim_{\nu\to
        \infty}\bigl[\omega\bigl(f(z),f(w_\nu)\bigr)
        -\omega\bigl(0,f(w_\nu)\bigr)\bigr]\cr
        &\le\liminf_{\nu\to\infty}\bigl[k_{\Delta^n}(z,w_\nu)
-\omega\bigl(0,f(w_\nu)
        \bigr)\bigr]\cr
        &\le\liminf_{\nu\to\infty}[k_{\Delta^n}(z,w_\nu)-k_{\Delta^n}
(0,w_\nu)]+
        \lim_{\nu\to\infty}\bigl[k_{\Delta^n}(0,w_\nu)-\omega\bigl(0,f(w_\nu)
        \bigr)\bigr]\cr
        &\le\limsup_{w\to x}[k_{\Delta^n}(z,w)-k_{\Delta^n}(0,w)]
+\mezzo\log\alpha\cr
        &<\mezzo\log(\alpha R),\cr}$$
that is~$f(z)\in E(\tau,\alpha R)$.

Finally, let $\sigma$ be a peculiar $x$-curve. 
Then \eqbJ\ implies that $f\bigl(
\sigma(t)\bigr)\in E(\tau,R)$ eventually for all~$R>0$, and this may happen
iff~$f\bigl(\sigma(t)\bigr)\to\tau$ as~$t\to1^-$, and we are done.\qedn

If the $\liminf$ in \eqliminf\ is equal to~$\mezzo\log\alpha$, we shall say
that~$f$
is {\sl $\alpha$-Julia} at~$x$.

It turns out that to compute \eqliminf\ it suffices to check what happens along
the image of~$\phe_x$:

\newthm Lemma \qu: Let $f\in\Hol(\Delta^n,\Delta)$, and take~$x\in\de\Delta^n$.
Then
$$\liminf_{w\to x}{1-|f(w)|\over1-\ln w\rn}=\liminf_{t\to1^-}{1-\bigl|f\bigl(
        \phe_x(t)\bigr)\bigr|\over1-t}.\neweq\eqqu$$

\pf Let us call~$\alpha$ the left-hand side of~\eqqu, and~$\beta$ the
right-hand side. Since~$\phe_x(t)\to x$ as~$t\to1^-$, clearly~$\alpha\le\beta$;
in particular, if~$\alpha=+\infty$ we are done. So assume~$\alpha<+\infty$; we
must show that~$\beta\le\alpha$.

Since $\alpha$ is finite, we can apply Theorem~\Julia\ (and its proof). So
there is~$\tau\in\de\Delta$ such that
$$\forevery{z\in\Delta^n}
{|\tau-f(z)|^2\over1-|f(z)|^2}\le\alpha\,\max_{|x_j|=1}
        \left\{{|x_j-z_j|^2\over1-|z_j|^2}\right\}.$$
Now, if~$\zeta\in\Delta$ then for $z=\phe_x(\zeta)=\zeta x$ one has
$$\max_{|x_j|=1}\left\{{|x_j-z_j|^2\over1-|z_j|^2}\right\}={|1-\zeta|^2\over
        1-|\zeta|^2};$$
therefore
$$\sup_{\zeta\in\Delta}\left\{{\bigl|\tau-f\bigl(\phe_x(\zeta)\bigr)\bigr|^2
        \over1-\bigl|f\bigl(\phe_x(\zeta)\bigr)\bigr|^2}\bigg/{|1-\zeta|^2\over
        1-|\zeta|^2}\right\}\le\alpha.\neweq\eqqub$$
Set $t_k=(k-1)/(k+1)$ for every~$k\in\N$. Clearly $t_k\in\Delta$
and~$t_k\to1^-$ as~$k\to+\infty$; moreover, $|1-t_k|^2/(1-|t_k|^2)=1/k$. It
follows that~$f\bigl(\phe_x(t_k)\bigr)\in\bar{E(\tau,\alpha/k)}$, by~\eqqub.
Now, $E(\tau,\alpha/k)$ is an euclidean disk of radius~$\alpha/(k+\alpha)$;
therefore
$$1-\bigl|f\bigl(\phe_x(t_k)\bigr)\bigr|\le
\bigl|\tau-f\bigl(\phe_x(t_k)\bigr)\bigr|\le{2\alpha\over k+\alpha}.$$
Since $1-|t_k|=2/(k+1)$ it follows that
$$\beta\le\limsup_{k\to\infty}{1-\bigl|f(\phe_x(t_k)\bigr)\bigr|\over1-|t_k|}
        \le\lim_{k\to\infty}\alpha\,{k+1\over k+\alpha}=\alpha,$$
and we are done.\qedn

In particular, thus, to check whether a given bounded holomorphic function~$f$
is \hbox{$\alpha$-Julia} at~$x\in\de\Delta^n$ it suffices to study the
function~$t\mapsto f(tx)$.
\smallbreak
{\it Remark 3.1:} The $\liminf$ \eqqu\ is always positive. Indeed,
$$\omega\bigl(0,f(w)\bigr)\le\omega\bigl(0,f(0)\bigr)+\omega\bigl(f(0),f(w)
        \bigr)\le\omega\bigl(0,f(0)\bigr)+k_{\Delta^n}(0,w);$$
therefore
$$k_{\Delta^n}(0,w)-\omega\bigl(0,f(w)\bigr)\ge-\omega\bigl(0,f(0)\bigr)>
        -\infty$$
and
$${1-|f(w)|\over1-\ln w\rn}\ge{1-|f(0)|\over2(1+|f(0)|)}>0$$
for all $w\in\Delta^n$.

\smallsect 4. The Julia-Wolff-Carath\'eodory theorem

We are finally ready to state and prove the Julia-Wolff-Carath\'eodory
Theorem~\zJWCA\ for bounded holomorphic functions in the
polydisk. 
If $v\in\C^n$
and~$f\colon\Delta^n\to\C$, we set
$${\de f\over\de v}(z)=df_z(v)=\sum_jv^j{\de f\over\de z^j}(z).$$

\newthm Theorem \JWC: Let $f\in\Hol(\Delta^n,\Delta)$ be a bounded holomorphic
function, and~$x\in\de\Delta^n$. Assume there is~$\alpha>0$ such that
$$\liminf_{w\to x}{1-|f(w)|\over1-\ln w\rn}=\alpha<+\infty.$$
Let $\tau\in\de\Delta$ be the restricted $E$-limit of~$f$ at~$x$, as given by
Theorem~\Julia. Then:
{$$\rKlim_{z\to x}{\tau-f(z)\over1-\tilde p_x(z)}=\alpha\tau;
        \leqno{\,\,\rm(i)}$$
\itm{(ii)}If $x_j$ is a \v Silov component of~$x$, then
$$\rKlim_{z\to x}{\tau-f(z)\over x_j-z_j}=\alpha\tau\bar{x_j};$$
$$\rKlim_{z\to x}{\de f\over\de x}(z)=\rKlim_{z\to x}{\de f\over\de\check x}(z)
        =\alpha\tau;\leqno{\!\rm(iii)}$$
\itm{(iv)}If $x_j$ is an internal component of~$x$, then
$$\rKlim_{z\to x}{\de f\over\de z_j}(z)=0;$$
\itm{(v)}If $x_j$ is a \v Silov component of~$x$, then $\de f/\de z_j$ has
restricted $K$-limit at~$x$.}

\smallbreak
{\it Remark 4.1:} If $x_j$ is an internal component of~$x$, that is~$|x_j|<1$,
then the incremental ratio~$\bigl(\tau-f(z)\bigr)/(x_j-z_j)$ is {\it not}
well-defined, because there are~$z\in\Delta^n$ with~$z_j=x_j$.
\smallbreak
{\it Remark 4.2:} If $x_j$ is a \v Silov component of~$x$, then $\de f/\de z_j$
might have a restricted $K$-limit at~$x$ different from the restricted
$K$-limit of the corresponding incremental ratio. For instance, choose $0<\beta
<\alpha<1$ and let $f\in\Hol(\Delta^2,\Delta)$ be given by
$$f(z_1,z_2)=1+\mezzo(\alpha+\beta)(z_1-1)+\mezzo(\alpha-\beta)(z_2-1).$$
Then it is easy to check that $f$ is $\alpha$-Julia at~$x=(1,1)$ but
$${\de f\over\de z_1}\equiv{\alpha+\beta\over2}$$
is different from~$\alpha$. Notice that, on the other hand, $\de f/\de x
\equiv\alpha$, as it should be.
\smallbreak
{\it Remark 4.3:} In general, as we shall see in Propositions~4.8,~4.9 and
Remark~4.6, the restricted $K$-limit of~$\de f/\de z_j$ at~$x$ is of the
form~$\beta_j\tau\bar{x_j}$, where all $\beta_j$'s are
non-negative, 
$\beta_j=0$
if $x_j$ is an internal component of~$x$, and~$\beta_1+\cdots+\beta_n=\alpha$.
\smallbreak

The proof of Theorem~\JWC\ will fill the rest of this section. The idea is to
show that the given functions are $K$-bounded, have limit along a special
restricted curve --- usually $t\mapsto\phe_x(t)$ --- and then apply
Theorem~\LinK. We begin with:

\newthm Lemma \qt: Let $f\in\Hol(\Delta^n,\Delta)$ be $\alpha$-Julia at~$x\in
\de\Delta^n$, and let~$\tau\in\de\Delta$ be its restricted $E$-limit at~$x$.
Fix~$M>1$. Then for all $z\in H(x,M)$ we have
$$\left|{\tau-f(z)\over1-\tilde p_x(z)}\right|\le2\alpha M^2\qquad\hbox{and}
        \qquad\left|{\tau-f(z)\over x_j-z_j}\right|\le2\alpha M^2,$$
where $x_j$ is any \v Silov component of~$x$.

\pf Take $z\in H(x,M)$ and set
$$\mezzo\log R=\log M-k_{\Delta^n}(0,z).\neweq\eqqt$$
Clearly, $z\in E(x,R)$; therefore $f(z)\in E(\tau,\alpha R)$ and thus
$$\lim_{s\to1^-}\bigl[\omega\bigl(f(z),s\tau\bigr)-\omega(0,s\tau)\bigr]
        -\omega\bigl(0,f(z)\bigr)<\log(\alpha R)$$
(notice that $-\omega\bigl(0,f(z)\bigr)\le\omega(f(z),s\tau)-\omega(0,s\tau)$
for all~$s<1$). So
$$\log(\alpha R)>\mezzo\log{|\tau-f(z)|^2\over1-|f(z)|^2}-\mezzo\log{1+|f(z)|
        \over1-|f(z)|}=\log{|\tau-f(z)|\over1+|f(z)|},$$
and \eqqt\ yields
$$\log{|\tau-f(z)|\over1+|f(z)|}<\log\alpha
+\log M^2-\log{1+\ln z\rn\over1-\ln z
        \rn}\le\log(\alpha M^2)-
\log{1+|\tilde p_x(z)|\over1-|\tilde p_x(z)|},$$
because $|\tilde p_x(z)|\le\ln z\rn$. Hence
$$\log\left|{\tau-f(z)\over1-\tilde p_x(z)}\right|\le\log{|\tau-f(z)|\over
        1-|\tilde p_x(z)|}
<\log\left(\alpha M^2{1+|f(z)|\over1+|\tilde p_x(z)|}\right)
        \le\log(2\alpha M^2).$$
Analogously, if $x_j$ is a \v Silov component of~$x$ we get
$$\log\left|{\tau-f(z)\over x_j-z_j}\right|\le\log{|\tau-f(z)|\over
        1-\ln z\rn}<\log\left(\alpha M^2{1+|f(z)|\over1+\ln z\rn}\right)
        \le\log(2\alpha M^2).$$
\qedn

\newthm Lemma \qd: Let $f\in\Hol(\Delta^n,\Delta)$ be $\alpha$-Julia at~$x\in
\de\Delta^n$, and let~$\tau\in\de\Delta$ be its restricted $E$-limit
at~$x$. Then
$$\lim_{t\to1^-}{\tau-f\bigl(\phe_x(t)\bigr)\over1-t}=\lim_{t\to1^-}(f\circ
        \phe_x)'(t)=\alpha\tau.$$

\pf Indeed, Lemma~\qu\ shows that $f\circ\phe_x$ is $\alpha$-Julia at~$1\in\de
\Delta$, and the assertion follows from the classical
Julia-Wolff-Carath\'eodory Theorem~\zJWC.\qedn

And so:

\newthm Corollary \qq: Let $f\in\Hol(\Delta^n,\Delta)$ be $\alpha$-Julia
at~$x\in
\de\Delta^n$, and let~$\tau\in\de\Delta$ be its restricted $E$-limit
at~$x$. Then
$$\rKlim_{z\to x}{\tau-f(z)\over1-\tilde p_x(z)}=\alpha\tau\qquad\hbox{and}
        \qquad\rKlim_{z\to x}{\tau-f(z)\over x_j-z_j}=\alpha\tau\bar{x_j},$$
where $x_j$ is any \v Silov component of~$x$.

\pf The first limit follows immediately from Lemmas~\qt, \qd\ and
Theorem~\LinK. For the second limit it suffices to remark that
$${\tau-f(z)\over x_j-z_j}={\tau-f(z)\over 1-\tilde p_x(z)}\,{1-\tilde p_x(z)
        \over x_j-z_j}$$
and
$${1-\tilde p_x\bigl(\phe_x(t)\bigr)\over
        x_j-\bigl(\phe_x(t)\bigr)_j}=\bar{x_j},$$
and then apply again Lemmas~\qt, \qd\ and Theorem~\LinK.\qedn

We have proved parts (i) and (ii) of Theorem~\JWC; furthermore, by
Lemma~\qd\ we know that $\de f/\de x$ has limit~$\alpha\tau$ along the
$x$-curve $t\mapsto\phe_x(t)$. So we must now deal with the $K$-boundedness of
the partial derivatives. To do so we need two further lemmas:

\newthm Lemma \qs: Take $M_1>M>1$ and set $r=(M_1-M)/(M_1+M)<1$. Fix $x\in\de
\Delta^n$ and let $\psi\in\Hol(\Delta,\Delta^n)$ be a complex geodesic such
that $z_0=\psi(0)\in H(x,M)$. Then $\psi(\Delta_r)\subset H(x,M_1)$.

\pf Let $\delta=\mezzo\log(M_1/M)>0$; then $\zeta\in\Delta_r$ iff $\omega(0,
\zeta)<\delta$. Then
$$\eqalign{\lim_{s\to1^-}\bigl[k_{\Delta^n}\bigl(\psi(\zeta),\phe_x(s)\bigr)&-
        \omega(0,t)\bigr]+k_{\Delta^n}\bigl(0,\psi(\zeta)\bigr)\cr
        &\le2k_{\Delta^n}
        \bigl(z_0,\psi(\zeta)\bigr)
+\lim_{s\to1^-}\bigl[k_{\Delta^n}\bigl(z_0,\phe_x(s)
        \bigr)-\omega(0,t)\bigr]+k_{\Delta^n}(0,z_0)\cr
        &<2\omega(0,\zeta)+\log M<\log M_1\cr}$$
for all $\zeta\in\Delta_r$.\qedn

We shall denote by $\kappa_{\Delta^n}\colon\Delta^n\times\C^n\to\R^+$ the
Kobayashi metric of~$\Delta^n$; it is well-known (see, e.g.,
[J-P,~Example~3.5.6]) that
$$\kappa_{\Delta^n}(z;v)=\max_j\left\{{|v_j|\over1-|z_j|^2}\right\}.$$

\newthm Lemma \km: Take $x\in\de\Delta^n$ and $M>1$. 
Then for any $z\in H(x,M)$
and $v\in\C^n$ we have
$$|1-\tilde p_x(z)|\kappa_{\Delta^n}(z;v)\le 2M^3\ln v\rn.\neweq\eqkm$$

\pf First of all,
$$(1-\ln z\rn)\kappa_{\Delta^n}(z;v)=\max_j\left\{{|v_j| (1-\ln
        z\rn)\over1-|z_j|^2}\right\}\le\ln v\rn;$$
moreover if $z\in H(x,M)$ we have
$${|1-\tilde p_x(z)|\over1-\ln z\rn}\le M\,{1-\ln p_x(z)\rn\over1-\ln z\rn},$$
by Corollary~\add. Now, being $z\in H(x,M)$, for any $\eps>0$ there is~$t<1$
such that
$$\eqalign{\eps+\log M&>k_{\Delta^n}\bigl(z,\phe_x(t)\bigr)-\omega(0,t)+
        k_{\Delta^n}(0,z)\cr
        &\ge k_{\Delta^n}\bigl(p_x(z),\phe_x(t)\bigr)-k_{\Delta^n}\bigl(
        \phe_x(t),0\bigr)+k_{\Delta^n}(0,z)\ge
        k_{\Delta^n}(0,z)-k_{\Delta^n}\bigl(0,p_x(z)\bigr);\cr}$$
therefore
$$\log M\ge k_{\Delta^n}(0,z)-k_{\Delta^n}\bigl(0,p_x(z)\bigr)\ge
        \mezzo\log{1-\ln p_x(z)\rn\over2(1-\ln z\rn)},$$
and \eqkm\ follows.\qedn
        
Then:

\newthm Proposition \qst: Let $f\in\Hol(\Delta^n,\Delta)$ be $\alpha$-Julia
at~$x\in
\de\Delta^n$, and let~$\tau\in\de\Delta$ be its restricted $E$-limit
at~$x$. Then
for every~$v\in\C^n$ the partial derivative $\de f/\de v$ 
is $K$-bounded at~$x$.

\pf Fix $v\in\C^n$, with $v\ne0$, and for every $z\in\Delta^n$ let $\psi_z\in
\Hol(\Delta,\Delta^n)$ be a complex geodesic with~$\psi_z(0)=z$ 
and $\psi'_z(0)=
v/\kappa_{\Delta^n}(z;v)$. Clearly, we can choose~$\psi_z$ 
of the form~$\gamma
\circ\phe_y$ for suitable~$y\in\de\Delta^n$ and $\gamma$ automorphism
of~$\Delta^n$.

Choose $r\in(0,1)$. Cauchy's formula yields
$$\eqalign{{\de f\over\de v}(z)
&=\kappa_{\Delta^n}(z;v)(f\circ\psi_z)'(0)=      
{\kappa_{\Delta^n}(z;v)\over 2\pi i}\int_{|\zeta|=r}{f\bigl(\psi_z(\zeta)\bigr)
        \over\zeta^2}\,d\zeta\cr
        &={1\over2\pi}\int_0^{2\pi}{f\bigl(\psi_z(re^{1\theta})\bigr)-\tau\over
        \tilde p_x\bigl(\psi_z(re^{i\theta})\bigr)-1}\cdot{\tilde
        p_x\bigl(\psi_z(re^{i\theta})\bigr)-1\over\tilde
        p_x(z)-1}\cdot{\bigl(\tilde p_x(z)-1\bigr)\kappa_{\Delta^n}(z;v)\over r
        e^{i\theta}}\,d\theta,\cr}\neweq\eqint$$
where we used the fact that $\int_0^{2\pi}e^{-i\theta}\,d\theta=0$.

Fix $M>1$ and take $z\in H(x,M)$. Choose $M_1>M$ so that
$(M_1-M)/(M_1+M)>r$. By Lemma~\qs\ we have~$\psi_z(\bar{\Delta_r})\subset
H(x,M_1)$; so Lemma~\qt\ yields
$$\left|{f\bigl(\psi_z(re^{i\theta})\bigr)-\tau\over
        \tilde p_x\bigl(\psi_z(re^{i\theta})\bigr)-1}
\right|\le 2\alpha M^2_1,$$
and the first factor in \eqint\ is bounded. For the second factor we first
remark that
$$\left|{\tilde p_x\bigl(\psi_z(re^{i\theta})\bigr)-1\over\tilde p_x(z)-1}
        \right|\le{\bigl|1-\tilde
        p_x\bigl(\psi_z(re^{i\theta})\bigr)\bigr|\over1-\bigl|
        \tilde p_x\bigl(\psi_z(re^{i\theta})\bigr)\bigr|}\cdot{1-\bigl|
   \tilde p_x\bigl(\psi_z(re^{i\theta})\bigr)\bigr|\over1-|\tilde p_x(z)|}\le
        M_1{1-\ln p_x\bigl(\psi_z(re^{i\theta})\bigr)\rn\over1-\ln p_x(z)
        \rn},$$
by Corollary~\add\ and Lemma~\qs. Now
$$\eqalign{\mezzo\log{1-\ln p_x\bigl(\psi_z(re^{i\theta})\bigr)\rn\over2(1-\ln
        p_x(z)\rn)}&\le k_{\Delta^n}\bigl(0,p_x(z)\bigr)-k_{\Delta^n}\Bigl(0,
        p_x\bigl(\psi_z(re^{i\theta})\bigr)\Bigr)\cr
   &\le k_{\Delta^n}\Bigl(p_x(z),p_x\bigl(\psi_z(re^{i\theta})\bigr)\Bigr)\le
k_{\Delta^n}\bigl(z,\psi_z(re^{i\theta})\bigr)=\mezzo\log{1+r\over1-r},\cr}$$
and so the second factor in \eqint\ is bounded too. By Lemma~\km, the third
factor is bounded by~$2M^3\ln v\rn/r$, and we are done.\qedn

{\it Remark 4.4:} Putting all together we have actually proved that
$$\left|{\de f\over\de v}(z)\right|\le C\alpha M^6\ln v\rn$$
for all $v\in\C^n$ and $z\in H(x,M)$, where $C$ is a universal constant,
obtained choosing the best~$r\in(0,1)$.
\smallbreak
Proposition~\qst\ together with
Lemma~\qd\ yield
$$\rKlim_{z\to x}{\de f\over\de x}(z)=\alpha\tau;$$
therefore to end the proof of Theorem~\JWC\ it suffices to show that parts~(iv)
and~(v) hold. This is accomplished in the following propositions.

\newthm Proposition \qo: Let $f\in\Hol(\Delta^n,\Delta)$ be $\alpha$-Julia
at~$x\in\de\Delta^n$, and let~$\tau\in\de\Delta$ be its restricted $E$-limit
at~$x$. Take $v\in\C^n$ with no \v Silov components with respect to~$x$. Then
$$\rKlim_{z\to x}{\de f\over\de v}(z)=0.$$

\pf For $t\in(0,1)$ let $\psi_t\in\Hol(\C,\C^n)$ be given by
$$\psi_t(\zeta)=tx+\zeta v.$$
We have
$$\ln\psi_t(\zeta)\rn=\max_j\{|tx_j+\zeta v_j|\}\le\max\{t,t\ln\interior x\rn
        +|\zeta|\ln v\rn\};$$
therefore if $r_t=(1-t)/\ln v\rn$ we have
$\psi_t(\bar{\Delta_r})\subset\Delta^n$.

For any $\theta\in\R$ put $\sigma^\theta(t)=\psi_t(r_te^{i\theta})$. Clearly,
$\sigma^\theta$ is a $x$-curve. Furthermore, $\sigma^\theta_x(t)=tx$, because
$v$ has no \v Silov components with respect to~$x$; so
$\sigma^\theta$ is restricted and (cf. Proposition~\Jafari) special.

Now, $\psi_t(0)=tx=\phe_x(t)$ and $\psi'_t(0)=v$ for any $t\in(0,1)$. Hence
$$\eqalign{{\de f\over\de v}(tx)=(f\circ\psi_t)'(0)
&={1\over2\pi i}\int_{|\zeta|
 =r_t}{f\circ\psi_t(\zeta)\over\zeta^2}\,d\zeta={1\over2\pi}\int_0^{2\pi}
  {f\bigl(\psi_t(r_te^{i\theta})\bigr)-\tau\over r_te^{i\theta}}\,d\theta\cr
  &={1\over2\pi}\int_0^{2\pi}{f\bigl(\sigma^\theta(t)\bigr)-\tau\over\tilde p_x
  \bigl(\sigma^\theta(t)\bigr)-1}\cdot{\tilde p_x\bigl(\sigma^\theta(t)\bigr)-1
        \over r_te^{i\theta}}\,d\theta.\cr}$$
Since $\sigma^\theta$ is a special restricted $x$-curve, we know that the first
factor in the integrand converges boundedly to~$\alpha\tau$ as~$t\to1^-$. For
the second factor,
$${\tilde p_x\bigl(\sigma^\theta(t)\bigr)-1
\over r_te^{i\theta}}=-{\ln v\rn\over
        e^{i\theta}};$$
therefore
$$\lim_{t\to 1^-}{\de f\over\de v}(tx)=-{\alpha\tau\ln
        v\rn\over2\pi}\int_0^{2\pi}{d\theta\over e^{i\theta}}=0,$$
and the assertion follows from Proposition~\qst\ and Theorem~\LinK.\qedn

{\it Remark 4.5:} The curve $\sigma^\theta(t)=tx+r_te^{i\theta}v$, where
$r_t=(1-t)/\ln v\rn$, is special iff~$v=\lambda\check x+\interior{v}$, where
$\lambda\in\C$ and $\interior{v}$ has no \v Silov components with respect
to~$x$, whereas in the proof of~[J, Theorem~5.(d)] it is mistakenly assumed to
be special (that is, tangent to the diagonal) for all~$v\in\C^n$.
\smallbreak

In particular, then, the case of \v Silov components requires a completely
different proof:

\newthm Proposition \atlast: Let $f\in\Hol(\Delta^n,\Delta)$ be $\alpha$-Julia
at~$x\in\de\Delta^n$, and let~$\tau\in\de\Delta$ be its restricted $E$-limit
at~$x$. Take a \v Silov component~$x_j$ of~$x$. Then $\de f/\de z_j$ has
restricted $K$-limit~$\beta\tau\bar{x_j}$ at~$x$, where~$\beta\ge0$.

\pf Let $M$ be the set of all $k=(k_1,\ldots,k_n)\in\N^n$ with $k_1,\ldots,k_n$
relatively prime and~$|k|=k_1+\cdots+ k_n>0$. Since the function
$(\tau+f)/(\tau-f)$ has positive real part, the generalized Herglotz
representation formula proved in~[K-K] yields
$${\tau+f(z)\over\tau-f(z)}=\sum_{k\in
        M}\left[\int\limits_{(\de\Delta)^n}{w^k+z^k\over
        w^k-z^k}\,d\mu_k(w)+C_k\right],\neweq\eqalu$$
for suitable $C_k\in\C$ and positive Borel measures~$\mu_k$ on~$(\de\Delta)^n$,
where $z^k=z_1^{k_1}\cdots z_n^{k_n}$ and~$w^k=w_1^{k_1}\cdots w_n^{k_n}$; the
sum is absolutely converging.

Let $X_k=\{w\in(\de\Delta)^n\mid w^k=x^k\}$, and set $\beta_k=\mu_k(X_k)\ge0$
and $\mu^o_k=\mu_k-\mu_k|_{X_k}$, where $\mu_k|_{X_k}$ is the restriction
of~$\mu_k$ to~$X_k$ (i.e., $\mu_k|_{X_k}(E)=\mu_k(E\cap X_k)$ for every Borel
subset~$E$). Notice that $X_k=\void$ (and so~$\beta_k=0$) as soon as~$k_j>0$
for some internal component~$x_j$ of~$x$.

Using these notations \eqalu\ becomes
$${\tau+f(z)\over\tau-f(z)}=\sum_{k\in M}\left[\beta_k\,{x^k+z^k\over x^k-z^k}
        +\int\limits_{(\de\Delta)^n}{w^k+z^k\over
        w^k-z^k}\,d\mu^o_k(w)+C_k\right].\neweq\eqald$$
In particular, if $z=tx=\phe_x(t)$ we get
$${\tau+f(tx)\over\tau-f(tx)}=\sum_{k\in M}\left[\beta_k\,{1+t^{|k|}\over
  1-t^{|k|}}+\int\limits_{(\de\Delta)^n} {w^k+t^{|k|}x^k\over w^k-t^{|k|}x^k}
        \,d\mu^o_k(w)+C_k\right].\neweq\eqalt$$
Let us multiply both sides by~$(1-t)$, and then take the limit as~$t\to1^-$.
The left-hand side, by Corollary~\qq, tends to~$2/\alpha$ (notice
that~$\alpha\ne0$, by Remark~3.1). For the right-hand side, 
first of all we have
$${1+t^{|k|}\over1-t^{|k|}}\,(1-t)={1+t^{|k|}\over1+\cdots+t^{|k|-1}}\to
        {2\over|k|}.$$
Next, if $|x^k|<1$ it is clear that
$$(1-t)\int\limits_{(\de\Delta)^n} {w^k+t^{|k|}x^k\over w^k-t^{|k|}x^k}
        \,d\mu^o_k(w)\to0.\neweq\eqalq$$
Otherwise, since $\mu^o_k(X_k)=0$, for every $\eps>0$ there exists an open
neighborhood~$A_\eps$ of~$X_k$ in~$(\de\Delta)^n$ such that~$\mu^o_k(A_\eps)<
\eps$. Then
$$\eqalign{(1-t)\biggl|\int\limits_{(\de\Delta)^n}&{w^k+t^{|k|}x^k\over
        w^k-t^{|k|}x^k}\,d\mu^o_k(w)\biggr|\cr
        &\le(1-t)\left|\int\limits_{A_\eps} {w^k+t^{|k|}x^k\over
  w^k-t^{|k|}x^k}\,d\mu^o_k(w)+\int\limits_{(\de\Delta)^n\setminus A_\eps}
        {w^k+t^{|k|}x^k\over w^k-t^{|k|}x^k}\,d\mu^o_k(w)\right|\cr
  &\le2\,{1-t\over1-t^{|k|}}\eps
+(1-t)\left|\int\limits_{(\de\Delta)^n\setminus A_\eps}
        {w^k+t^{|k|}x^k\over    w^k-t^{|k|}x^k}\,d\mu^o_k(w)\right|\to
        {2\over|k|}\eps.\cr}$$
Since this happens for all $\eps>0$, 
it follows that \eqalq\ holds in this case
too. Summing up, we have found
$${1\over\alpha}=\sum_{k\in M}{\beta_k\over|k|};\neweq\eqals$$
in particular, the series on the right-hand side is converging.

Differentiating \eqald\ with respect to~$z_j$ we get
$${\de f\over\de
    z_j}(z)=\bar{\tau}\bigl(\tau-f(z)\bigr)^2\sum_{k\in M}
k_j{z^k\over z_j}\left[\beta_k{x^k\over(x^k-z^k)^2}+\int\limits_{(\de\Delta)^n}
        {w^k\over(w^k-z^k)^2}\,d\mu^o_k(w)\right].\neweq\eqalst$$
Since we
already know that $\de f/\de z_j$ is $K$-bounded, 
it suffices to show that $\de
f/\de z_j$ has limit along the $x$-curve $t\mapsto tx$, where we have
$$\eqalign{{\de f\over\de
        z_j}(tx)=\bar\tau\left({\tau-f(tx)\over1-t}\right)^2\sum_{k\in M}k_j
        t^{|k|-1}\bar{x_j}\biggl[&\beta_k\left({1-t\over1-t^{|k|}}\right)^2\cr
        &+x^k(1-t)^2\int\limits_{(\de\Delta)^n}{w^k\over(w^k-t^{|k|}x^k)^2}
        \,d\mu^o_k(w)\biggr].\cr}$$
The same argument used before shows that
$$(1-t)^2\int\limits_{(\de\Delta)^n}{w^k\over(w^k-t^{|k|}x^k)^2}\,d\mu^o_k(w
)\to 0$$
as $t\to1^-$. Therefore
$$\lim_{t\to1^-}{\de f\over\de z_j}(tx)=\alpha^2\tau\bar{x_j}\sum_{k\in M}
        \beta_k{k_j\over|k|^2},\neweq\eqalo$$
where the series is converging because $\beta_k k_j/|k|^2\le\beta_k/|k|$, and
we are done.\qedn

{\it Remark 4.6:} If $x_j$ is an internal component of~$x$, then the sum in
\eqalo\ vanishes. Indeed, we have already remarked that $\beta_k=0$ if~$k_j>0$,
and so in this case~$\beta_k k_j=0$ always. Furthermore, \eqalo\ and~\eqals\
yield
$$\lim_{t\to1^-}{\de f\over\de x}(tx)=\sum_{j=1}^nx_j\lim_{t\to1^-}{\de f\over
  \de z_j}(tx)=\alpha^2\tau\sum_{j=1}^n|x_j|^2\sum_{k\in M}\beta_k{k_j\over
        |k|^2}=\alpha^2\tau\sum_{k\in M}{\beta_k\over|k|}=\alpha\tau,$$
(where we used again the fact that $\beta_k k_j=0$ always if $|x_j|<1$), as it
should be according to Theorem~\JWC.(iii) and Remark~4.3.
\smallbreak

{\it Remark 4.7:} If $d_x=1$, then $\de f/\de v$ has restricted
$K$-limit at~$x$ for all~$v\in\C^n$. Indeed, in this case all $v\in\C^n$ are of
the form~$\lambda\check x+\interior{v}$, where $\interior{v}$ has no \v
Silov components with respect to~$x$, and the assertion follows from
Theorem~\JWC.(iii) and~(iv).

\smallsect 5. The multidimensional case

As discussed in the introduction, given the correct setup and enough
geometrical information, it is possible to obtain
Julia-Wolff-Carath\'eodory-like theorems for holomorphic maps between any
kind of
domains. In~[A3,~5] we discussed the situation for maps between strongly convex
and strongly pseudoconvex domains; in the previous section we studied the
situation for functions from a polydisk into the unit disk in~$\C$. This
section is devoted to describe what happens for maps from a polydisk into
another polydisk, or for maps from a polydisk into a strongly (pseudo)convex
domain.

Let us start with a $f\in\Hol(\Delta^n,\Delta^m)$. The Julia condition
$$\liminf_{w\to x}\bigl[k_{\Delta^n}(0,w)-k_{\Delta^n}\bigl(0,f(w)\bigr)\bigr]
        =\mezzo\log\alpha<+\infty\neweq\eqJzero$$
translates in
$$\min_j\liminf_{w\to x}{1-|f_j(w)|\over1-\ln
w\rn}<+\infty.\neweq\eqJuno$$
In other words, \eqJzero\ is equivalent to assuming that at least one component
of~$f$ is \hbox{$\alpha$-Julia} for some~$\alpha$, without saying anything else
on the other components. Therefore if we write~$f=\check f+\interior f$, where
$\check f$ contains the components of~$f$ 
satisfying a Julia condition (possibly
with different~$\alpha$'s), and $\interior f$ the other components, we recover
for~$\check f$ results exactly like the one described in
Theorems~\Julia\ and~\JWC, whereas we cannot say anything about~$\interior f$.
For instance, if $g\in\Hol(\Delta,\Delta)$ is given by
$$g(\zeta)=\exp\left(-{\pi\over2}-i\log(1-\zeta)\right),$$
then~$g(t)$ has no limit as~$t\to1^-$, and the map
$$f(z_1,z_2)=\bigl(z_1,\mezzo g(z_2)\bigr)$$
satisfies \eqJuno\ with~$x=(1,1)$ and~$\alpha=1$, but $f_2(t,t)$ has no limit
as~$t\to1^-$. In particular, it is easy to check that
$$f\Bigl(E\bigl((1,1),1\bigr)\Bigr)\not\subseteq E\bigl((1,1),1\bigr),$$
and so in a general Julia lemma (like Theorem~\zJA) one is forced to consider
both small and big horospheres.
\smallbreak
{\it Remark 5.1:} Even if it is less natural 
from the point of view of geometric
function theory, one might of course
consider maps $f\in\Hol(\Delta^n,\Delta^m)$ satisfying
$$\max_j\liminf_{w\to x}{1-|f_j(w)|\over1-\ln w\rn}<+\infty$$
instead of \eqJuno. Then all components of~$f$ are $\alpha$-Julia (for possibly
different $\alpha$'s), and we recover Theorems~\Julia\ and~\JWC\ for all
components of~$f$.
\smallbreak
The situation is more interesting if we consider maps~$f\in\Hol(\Delta^n,D)$,
where~$D\subset\subset\C^m$ is a strongly convex domain. To state the result in
this case we need some preparation. Fix once for all a point~$z_0\in D$; then
for every~$y\in\de D$ there exists (see~[L, A2]) a unique complex
geodesic~$\psi_y\in\Hol(\Delta,D)\cap C^1(\bar{\Delta},\bar{D})$ such
that~$\psi_y(0)=z_0$ and~$\psi_y(1)=y$. Associated to~$\psi_y$ there is a
holomorphic retraction~$q_y\colon\bar D\to\psi_y(\bar\Delta)$ such that~$q_y
\circ q_y=q_y$; setting~$\tilde q_y=\psi^{-1}_y\circ q_y\colon\bar D\to\bar
\Delta$ we clearly have~$\tilde q_y\circ\psi_y=\id_\Delta$. In particular,
$\tilde q_y(y)=1$. Given~$f\in\Hol(\Delta^m,D)$ and~$y\in\de D$ we can
associate to~$f$ the function~$\tilde f_y=\tilde q_y\circ f$;
roughly speaking, $\tilde f_y$ is the component of~$f$ in the direction of~$y$,
and~$f-f_y$ (where~$f_y=q_y\circ f$) contains the components of~$f$ ``normal''
to~$y$.

Then the Julia-Wolff-Carath\'eodory theorem for this case is:

\newthm Theorem \JWCmix: Let $D\subset\subset\C^m$ be a bounded strongly
convex~$C^3$ domain, and fix~$z_0\in D$. Let~$x\in\de D$
and~$f\in\Hol(\Delta^n,D)$ be such that
$$\liminf_{w\to x}\bigl[k_{\Delta^n}(0,w)-k_D\bigl(z_0,f(w)\bigr)\bigr]
        =\mezzo\log\alpha<+\infty,\neweq\eqJF$$
where~$k_D$ is the Kobayashi distance of~$D$.
Then $f$ has restricted $E$-limit~$y\in\de D$ at~$x$, and the following maps
are $K$-bounded at~$x$:
{\smallskip
\itm{(i)}$[1-\tilde f_y(z)]/[1-\tilde p_x(z)]$;
\itm{(ii)}$[f(z)-f_y(z)]/[1-\tilde p_x(z)]^{1/2}$;
\itm{(iii)}$[1-\tilde f_y(z)]/[x_j-z_j]$, where $x_j$ is a \v Silov component
of~$x$;
\itm{(iv)}$[f(z)-f_y(z)]/[x_j-z_j]^{1/2}$, where $x_j$ is a \v Silov component
of~$x$;
\itm{(v)}$\de\tilde f_y/\de v$, for all~$v\in\C^n$;
\itm{(vi)}$[1-\tilde p_x(z)]^{1/2}d(f-f_y)_z(v)$, again for all~$v\in\C^n$.
\smallskip
\noindent Furthermore, function {\rm (i)} and function {\rm (v)}
for~$v=x$,~$\check x$ have restricted $K$-limit~$\alpha$ at~$x$; maps~{\rm
(ii), (iv)} and function~{\rm (v)} for $v$ without \v Silov components with
respect to~$x$ have restricted $K$-limit~$0$ at~$x$; function~{\rm (iii)} has
restricted $K$-limit~$\alpha\bar{x_j}$ at~$x$; function~{\rm (v)} for~$v=x_j$,
where $x_j$ is a \v Silov component of~$x$, 
has restricted $K$-limit at~$x$; and
map~{\rm (vi)} for~$v=\lambda\check x+\interior v$ (where
$\lambda\in\C$ and~$\interior v$ has no \v Silov 
components with respect to~$x$)
has restricted
$K$-limit~$0$ at~$x$.}

\pf The arguments needed are a mixture of the ones used to prove Theorem~\JWC\
and the ones used to prove the Julia-Wolff-Carath\'eodory theorem for strongly
convex domains (see~[A3]); so we shall just sketch the necessary modifications.

The existence of the restricted $E$-limit~$y\in\de D$ at~$x$ is proved exactly
as in Theorem~\Julia\ (cf.~[A3, Proposition~1.19]), and it follows from the
fact that horospheres in strongly convex domains touch the boundary in exactly
one point, as in the disk (see~[A1]).

Next, \eqJF\ implies
$$\eqalign{\liminf_{t\to1^-}\mezzo\log{1-|\tilde
   f_y(tx)|\over1-t}&=\liminf_{t\to1^-}\bigl[k_{\Delta^n}(0,tx)-k_D\bigl(z_0,
        f_y(tx)\bigr)\bigr]\cr
        &\ge\liminf_{t\to1^-}\bigl[k_{\Delta^n}(0,tx)-k_D\bigl(z_0,
        f(tx)\bigr)\bigr]\ge\mezzo\log\alpha.\cr}$$
On the other hand, since
$$\tilde f_y\circ\phe_x\bigl(E(1,R)\bigr)=\tilde p_y\Bigl(f\circ\phe_x\bigl(
        E(1,R)\bigr)\Bigr)
\subseteq\tilde p_y\bigl(E(y,\alpha R)\bigr)=E(1,\alpha R),$$
it follows that $\tilde f_y\circ\phe_x$ is $\alpha$-Julia; therefore
the classical Julia-Wolff-Carath\'eodory theorem implies
$$\limsup_{t\to 1^-}\mezzo\log{1-|\tilde f_y(tx)|\over1-t}\le\limsup_{t\to1^-}
        \mezzo\log{|1-\tilde f_y(tx)|\over1-t}\le\mezzo\log\alpha,$$
and $\tilde f_y\in\Hol(\Delta^n,\Delta)$ is $\alpha$-Julia at~$x$, by
Lemma~\qu. Hence we can apply Theorem~\JWC\ to~$\tilde f_y$, and all
the assertions concerning functions~(i), (iii) and~(v) are proved.

So we are left with $f-f_y$. The proof of~[A3, Proposition~3.7] applies word by
word, replacing~[A3, Proposition~3.4] by Lemma~\qt, and so maps~(ii) and~(iv)
are $K$-bounded at~$x$. Repeating the proof of~[A3, Proposition~3.9] we get the
restricted $K$-limit for map~(ii), and arguing as in the proof of
Corollary~\qq\ we deal with map~(iv) too.

Since, by Lemma~\km, $\bigl(1-\tilde p_x(z)\bigr)\kappa_{\Delta^n}(z;v)$
is $K$-bounded for all~$v\in\C^n$, the arguments in the proof
of~[A3,~Proposition~3.14] and Proposition~\qst\ show that the map~(vi) is
$K$-bounded for all~$v\in\C^n$. Finally, the proof of Proposition~\qo, suitably
adapted as in~[A3, Proposition~3.18] and recalling Remark~4.5, shows that the
map~(vi) has restricted \hbox{$K$-limit~0} at~$x$ for the indicated~$v$'s, as
claimed.\qedn

We end this paper by remarking that, using the techniques described in~[A5],
one can localize the statement of 
Theorem~\JWCmix\ near~$y\in\de D$, obtaining a
similar result for maps with values 
in a strongly pseudoconvex domain. We leave
the details to the interested reader.

\setref{Ju2}
\beginsection References

\art A1 Abate, M.: Horospheres and iterates of holomorphic maps! Math. Z.! 198
1988 225-238

\art A2 Abate, M.: Common fixed points of commuting holomorphic maps! Math.
Ann.! 283 1989 645-655

\art A3 Abate, M.: The Lindel\"of principle and the angular derivative in
strongly convex domains! J. Analyse Math.! 54 1990 189-228

\book A4 Abate, M.: Iteration theory of holomorphic maps on taut manifolds!
Rende: Me\-di\-ter\-ranean Press, 1989

\art A5 Abate, M.: Angular derivatives in strongly pseudoconvex domains! Proc.
Symp.  Pure Math.! {52, \rm Part 2,} 1991 23-40

\book B Burckel, R.B.: An introduction to classical complex analysis! New
York: Academic Press, 1979

\art C Carath\'eodory, C.: \"Uber die Winkelderivierten von beschr\"ankten
analytischen Funktionen! Sitzungsber. Preuss. Akad. Wiss., Berlin! {} 1929 39-54

\art C-K Cima, J.A., Krantz, S.G. : The Lindel\"of principle and normal
functions of several complex variables! Duke Math. J.! 50 1983 303-328

\art {\v C} \v Cirka, E.M.: The Lindel\"of and Fatou theorems in $\C^n$! Math.
USSR-Sb.! 21 1973 619-641

\art D Dovbush, P.V.: Existence of admissible limits of functions of several
complex variables! Sib. Math. J.! 28 1987 83-92

\art D-Z Dro\v z\v zinov, Ju.N., Zav'jalov, B.I.: On a multi-dimensional
analogue of a theorem of Lindel\"of! Sov. Math. Dokl.! 25 1982 51-52

\art H Herv\'e, M.: Quelques propri\'et\'es des applications analytiques d'une
boule \`a $m$ dimensions dans elle-m\^eme! J. Math. Pures Appl.! 42 1963 117-147

\art J Jafari, F. : Angular derivatives in polydiscs! Indian J. Math.! 35 1993
197-212

\book J-P  Jarnicki, M., Pflug, P.: Invariant distances and metrics in complex
analysis! Berlin: Walter de Gruyter, 1993

\art Ju1 Julia, G.: M\'emoire sur l'it\'eration des fonctions rationnelles! J.
Math. Pures Appl.! 1 1918 47-245

\art Ju2 Julia, G.: Extension nouvelle d'un lemme de Schwarz! Acta Math.! 42
1920 349-355

\art Kh Khurumov, Yu.V.: On Lindel\"of's theorem in $\C^n$! Soviet Math. Dokl.!
28 1983 806-809

\book K1 Kobayashi, S.: Hyperbolic manifolds and holomorphic mappings! New
York: Dekker, 1970

\art K2 Kobayashi, S.: Intrinsic distances, measures, and geometric function
theory! Bull. Amer. Math. Soc.! 82 1976 357-416

\art Ko Kor\'anyi, A.: Harmonic functions on hermitian hyperbolic spaces! Trans.
Amer. Math. Soc.! 135 1969 507-516

\art K-S Kor\'anyi, A., Stein, E.M.: Fatou's theorem for generalized
half-planes! Ann. Scuola Norm. Sup. Pisa! 22 1968 107-112

\art K-K Kosbergenov, S., Kytmanov, A.M.: Generalizations of the Schwarz and
Riesz-Her\-glotz formulas in Reinhardt domains. (Russian)! Izv. Vyssh. Uchebn.
Zaved. Mat.! {\rm n.~10} 1984 60-63

\art Kr Krantz, S.G.: Invariant metrics and the boundary behavior of holomorphic
functions on domains in $\C^n$! J. Geom. Anal.! 1 1991 71-97

\art L-V Landau, E., Valiron, G.: A deduction from Schwarz's lemma! J. London
Math. Soc.! 4 1929 162-163

\art L Lempert, L.: La m\'etrique de Kobayashi et la repr\'esentation  des
domaines sur la boule! Bull. Soc. Math. France! 109 1981 427-474

\art Li Lindel\"of, E.: Sur un principe g\'en\'erale de l'analyse et ses
applications \`a la theorie de la repr\'esentation conforme! Acta Soc. Sci.
Fennicae! 46 1915 1-35

\art M Minialoff, A.: Sur une propri\'et\'e des transformations dans l'espace de
deux variables complexes! C.R. Acad. Sci. Paris! 200 1935 711-713

\art N Nevanlinna, R.: Remarques sur le lemme de Schwarz! C.R. Acad. Sci. Paris!
188 1929 1027-1029

\book R Rudin, W.: Function theory in the unit ball of $\C^n$! Berlin: Springer,
1980

\book S Stein, E.M.: The boundary behavior of holomorphic functions of several
complex variables! Princeton: Princeton University Press, 1972

\art W Wlodarczyk, K.: Julia's lemma and Wolff's theorem for $J^*$-algebras!
Proc. Amer. Math. Soc.! 99 1987 472-476

\art Wo Wolff, J.: Sur une g\'en\'eralisation d'un th\'eor\`eme de Schwarz! C.R.
Acad. Sci. Paris! 183 1926 500-502

\bye